\documentclass{article}
\usepackage{amsmath}
\usepackage{amsthm}
\usepackage{amssymb}
\usepackage{amsfonts}
\usepackage{ mathrsfs }
\usepackage{enumerate}
\usepackage[utf8]{inputenc}
\usepackage[english]{babel}
\usepackage[ top = 3 cm, bottom = 2 cm]{geometry}

\newcommand{\nc}{\newcommand}
\nc{\ep}{\varepsilon}
\nc{\n}[1]{\mathscr{#1}}
\nc{\eps}[1]{{#1}_{\varepsilon}}
\nc{\be}{\begin{equation}}
\nc{\ee}{\end{equation}}
\nc{\m}[1]{\mathcal{#1}}
\nc{\bb}[1]{\mathbb{#1}}
\theoremstyle{plain}
\newtheorem{theorem}{Theorem}
\newtheorem{lemma}{Lemma}
\theoremstyle{definition}

\newtheorem{corollary}{Corollary}

\begin{document}
\title{Optimal Control of the Multiphase Stefan Problem\footnote{This research is funded by NSF grant \#1359074} }
\author{Ugur G. Abdulla and Bruno Poggi}
\date{}
\maketitle
\begin{center} {\large\noindent \textsc{Department of Mathematics, Florida Institute of Technology, Melbourne, Florida 32901}}
\par \medskip\bigskip\end{center}
{\bf Abstract.} We consider the inverse multiphase Stefan problem, where information on the heat flux on the fixed boundary is missing and must 
be found along with the temperature and free boundaries. Optimal control framework is pursued, where boundary heat flux is the control, and the optimality criteria consist of the minimization of the $L_2$-norm declination of the trace of the solution to the Stefan problem from the temperature measurement on the fixed right boundary. The state vector solves multiphase Stefan problem in a weak formulation, which is equivalent to Neumann problem for the quasilinear parabolic PDE with discontinuous coefficient. Full discretization through finite differences is implemented and discrete optimal control problem is introduced. We prove well-posedness in a Sobolev space framework and convergence of discrete optimal control problems to the original problem both with respect to the cost functional and control. Along the way, the convergence of the method of finite differences for the weak solution of the multiphase Stefan problem is proved. The proof is based on achieving a uniform $L_{\infty}$ bound, and $W_2^{1,1}$-energy estimate for the discrete multiphase Stefan problem.



\newpage
\section{Description of Main Results}\label{description of results}
\subsection{Introduction and Motivation}\label{E:1:1}
Consider the general multi-phase Stefan problem (\cite{LSU}): find the temperature function and phase transition boundaries $\{u(x,t),\quad\xi_j(t),j=\overline{1,J}\}$ in $D=\{0<x<\ell,\quad 0<t\leq T\}$ satisfying the following conditions:
\be\label{pde}
\alpha(u)u_t-(k(u)u_x)_x=f(x,t),\qquad (x,t)\in D, \quad u(x,t)\neq u^j,  j=\overline{1,J},
\ee
\begin{align}
&u|_{x=\xi_j(t)}=u^j, &\qquad 0<t\leq T, \quad j=\overline{1,J},\label{phasetrans} \\[3mm]
&[u]|_{x=\xi_j(t)}=0, &\qquad 0<t\leq T,\quad  j=\overline{1,J}, \label{ucont} \\[3mm]
&[k(u)u_x]|_{x=\xi_j(t)} = \gamma_j\frac{d\xi_j}{dt}, &\qquad 0<t\leq T,\quad  j=\overline{1,J}, \label{stefancond} \\[3mm]
&u(x,0) = \phi(x), &\qquad 0\leq x \leq \ell, \label{initial} \\[3mm]
&k(u)u_x|_{x=0}=g(t),\qquad k(u)u_x|_{x=\ell} = p(t),&\qquad 0<t\leq T, \label{flux}
\end{align}
where $\alpha, k$ are positive $C^1$ functions on each segment $(-\infty,u^1],$ $[u^j,u^{j+1}], j=1,\ldots,J-1$ and $[u^{J},+\infty)$ with 1st type discontinuity at $u=u^j, ~j=\overline{1,J}$, where $u^1<u^2<\cdots<u^J$ are known values; each $\gamma_j, j=\overline{1,J}$ is a known positive number, and $[u]|_{x=\xi_j}$ is the jump of $u$ at $\xi_j$, defined as
\[
[u]|_{x=\xi_j} = u|_{x=\xi_j}^+ - u|_{x=\xi_j}^-
\]
where $ u|_{x=\xi_j}^+$ (or $u|_{x=\xi_j}^-$) is the limit value of $u$ at $(\xi_j(t),t)$ taken in the region $\{(x,t): u>u^j\}$ (or $\{(x,t): u<u^j\}$). We define
\begin{equation}\label{initialinterfaceset}
\xi_j(0)=\Big\{x\in[0,\ell]~|~\phi(x)=u^j\Big\},\qquad j=\overline{1,J},
\end{equation}
that is, the phase transition boundaries at the initial time as level sets of the given initial function $\phi$.  In the physical context, $f$ characterizes the density of the sources, $\phi$ is the initial temperature, $g$ and $p$ are the heat fluxes on the left and right fixed boundaries respectively, each $u^j$ represents a phase transition temperature, and (\ref{stefancond}) is the Stefan condition expressing the conservation law according to which the free boundaries are pushed by the jump of the heat flux from different phases.

The weak formulation of the multiphase Stefan problem, as well as existence and uniqueness of the weak solution to the multiphase Stefan problem, was first proved in \cite{Kamenomostskaya, Oleinik}. We refer to monographs \cite{LSU, Meyrmanov} for the extensive list of references.

Assume now that some of the data is not available, or involves some measurement error. For example, suppose that the heat flux, $g$, at the fixed boundary $x=0$ is not known and must be found along with the temperature $u$ and the phase transition boundaries $\xi_j$. As compensation for not knowing this function, we must have access to additional information, which may come, for instance, as a measurement of the temperature at the fixed boundary $x=\ell$:
\begin{equation}\label{extra}
u(\ell,t)=\nu(t),\qquad 0<t\leq T.
\end{equation}

\textbf{Inverse Multiphase Stefan Problem (IMSP).} Find the functions $u(x,t)$, $~\xi_j(t),~ j=\overline{1,J},$ and the boundary heat flux $g(t)$ satisfying (\ref{pde})-(\ref{flux}),(\ref{extra}).

Motivation for the IMSP arose from the modeling of bioengineering problems on the laser ablation of biological tissues through a multiphase Stefan problem (\ref{pde})-(\ref{flux}). Laser ablation creates three phases - solid (skin),
fluid (melted skin) and gas (evaporated skin). Free boundaries $\xi_1(t)$ and $\xi_2(t)$ are measuring ablation depth separating solid/fluid and fluid/air regions at the moment $t$. The temperature measurement on the fixed boundary $x=0$ has an error, which makes it impossible to measure flux accurately. On the other side, temperature $\nu$ on the bottom fixed boundary $x=\ell$ can be measured, and the boundary flux $g$ can be identified by solving the IMSP. Moreover, our variational approach allows us to regularize errors in the measurement $\nu$. Our suggested approach to solve the IMSP is quite robust and can be applied for the identification of various functions such as $f, \alpha, k$. Laser ablation of the skin creates a complex environment with three phases, and phase transition boundaries are not necessarily functions of time, but may form complicated sets. The major advantage of our approach in this paper is that the weak formulation of the multiphase Stefan problem doesn't impose any explicit restriction on the structure of the free boundary, which is defined as a level set of the weak solution. Another important motivation arises from the optimal control of the laser ablation process, in which one tries to achieve particular temperature distribution on the bottom boundary $x=\ell$ by controlling flux on the fixed surface boundary, which is essential to guarantee that the ablation depth is not more than required, and no healthy tissue is affected.

The IMSP is not well posed in the sense of Hadamard. That is, if the data is not sufficiently coordinated, there may be no solution. Even if it exists, it might be not unique, and most importantly, there is no continuous dependence of the solution on the data functions. 

We refer to a recent paper  \cite{Abdulla1} for a review of the literature on Inverse Stefan Problems. The one-phase inverse Stefan problem (ISP) was first mentioned in \cite{Cannon3}, in which an unknown heat flux is to be determined under the given free boundary. The variational approach for solving this ill-posed inverse Stefan problem was used in \cite{BudakVasileva1, BudakVasileva2}. The first result on the optimal control of the Stefan problem appeared in \cite{Vasilev}, in which an optimal temperature along the fixed boundary must be determined to guarantee that the solutions of the Stefan problem stay close to the measurements taken at the final time. In \cite{Vasilev}, the existence result was proved. In \cite{Yurii}, the Fr\'echet derivative was found, the convergence of the finite difference scheme was proved, and Tikhonov regularization was suggested. Later development of the inverse Stefan problem proceeded in these two directions: Inverse Stefan problems with given phase boundaries were considered in \cite{Alifanov,Bell,Budak,Cannon,Carasso,Ewing1,Ewing2,Hoffman,Sherman,Goldman}; optimal control of Stefan problems, or equivalently inverse problems with unknown phase boundaries were investigated in \cite{Baumeister,Fasano,Hoffman1,Hoffman2,Jochum2,Jochum1,Knabner,Lurye,Nochetto, Niezgodka,Primicero,Sagues,Talenti,Goldman}. We refer to the monography \cite{Goldman} for a complete list of references of both types of inverse Stefan problems, both for linear and quasilinear parabolic equations.

In two recent papers \cite{Abdulla1, Abdulla2}, a new variational formulation of the one-phase ISP was developed. Optimal control framework was implemented, in which the boundary heat flux and the free boundary are components of the control vector and optimality criteria consist of the minimization of the sum of $L_2$-norm declinations from the available measurement of the temperature on the fixed boundary and available 
information on the phase transition temperature on the free boundary. This approach allows 
one to tackle situations when the phase transition temperature is not known explicitly, and is available through measurement with possible error. It also allows for the development of iterative numerical methods of least computational cost due to the fact that for every given control vector, the parabolic PDE is solved in a fixed region instead of a full free boundary problem. In \cite{Abdulla1}, the well-posedness in a Sobolev space framework and 
convergence of time-discretized optimal control problems is proved. In \cite{Abdulla2}, full discretization is implemented and the convergence of the discrete optimal control problems to the original problem both with respect to cost functional and control is proved.  The main advantage of this method is that numerically at each step, the problem to be solved is only a Neumann problem, and not a full free boundary problem. Moreover, the Neumann condition replaces the Stefan condition on the free boundary. In recent papers \cite{Abdulla3, Abdulla4}, the Fr\'echet differentiability and first order optimality condition in Besov spaces framework is proved and the formula for the Fr\'echet gradient is derived. 

This approach is not applicable to multiphase Stefan problem. The reason is that the Stefan condition on the phase transition boundary includes the flux calculated from both phases. Therefore, it can't be treated as a Neumann condition, even if we include the free boundary as one of the control components. In the current paper we develop a new approach based on the weak formulation of the multiphase Stefan problem, as a boundary value problem for the nonlinear PDE with discontinuous coefficients in a fixed domain. The main goal of this paper is to solve the IMSP in the optimal control framework by employing the weak formulation of the multiphase Stefan problem. We prove the existence of the optimal control and convergence of the sequence of discrete optimal control problems to the continuous problem both with respect to functional and control. The proof is based on the proof of uniform $L_{\infty}$ bound, and $W_2^{1,1}$-energy estimate for the discrete multiphase Stefan problem. We address the problem of Fr\'echet differentiability and application of the iterative gradient methods in Hilbert spaces in an upcoming paper.

In Section~\ref{E:1:1a} we describe the notation of Sobolev spaces used in this paper. In Section~\ref{E:1:1b}  we formulate IMSP as an optimal control problem. In Section~\ref{E:1:1c} we perform full discretization through finite differences and formulate discrete optimal control problem. In Section~\ref{E:1:1d} the main results are formulated.  In Section~\ref{E:2:1} we prove the existence and uniqueness of the discrete state vector. We present the proof of the main results in Section~\ref{E:3:1}. In Section~\ref{E:3:1a} we prove $L_{\infty}$ estimation for the discrete multiphase Stefan problem. In Section~\ref{E:3:1b} we prove $W_2^{1,1}$-energy  estimation for the discrete multiphase Stefan problem. Based on these estimations we prove the existence of the optimal control in Section~\ref{E:3:1c}. Proof of the convergence of the discrete optimal control problems to continuous optimal control problem is completed in Section~\ref{E:3:1d}.

\subsection{Notation of Sobolev Spaces}\label{E:1:1a}

$L_2(0,T)$ - Space of Lebesgue square-integrable functions. It is a Hilbert space with inner product
\[ (u,v)=\int_0^T uv\,dt. \]
$L_{\infty}(0,T)$ - Space of essentially bounded functions. It is a Banach space with norm
\[
\Vert u\Vert_{L_{\infty}[0,T]} =\underset{0\leq t\leq T}{\text{ esssup }}|u(t)|.
\]
$W_2^k(0,T), k=1,2,...$ - Hilbert space of all elements of $L_2(0,T)$ whose weak derivatives up to order $k$  exist and belong to $L_2(0,T)$. The inner product is defined as
\[ (u,v)=\int_0^T \sum_{s=0}^k \frac{d^su}{dt^s} \frac{d^sv}{dt^s} \,dt.  \]
$L_2(D)$ - Hilbert space with inner product
\[ (u,v)=\int_{D} uv\, dx\, dt . \]
$W_2^{1,0}(D)$ - Hilbert space of all elements of $L_2(D)$ that have a weak derivative in the $x$ direction, $\frac{\partial u}{\partial x}$, and such that it belongs to $L_2(D)$. The inner product is defined as
\[ (u,v)=\int_{D} \Big ( uv + \frac{\partial u}{\partial x}\frac{\partial v}{\partial x} \Big ) \,dx\, dt. \]
$W_2^{1,1}(D)$ - Hilbert space of all elements of $L_2(\Omega)$ with weak derivatives of first order, $\frac{\partial u}{\partial x}$, $\frac{\partial u}{\partial t}$. Also its weak derivatives must belong to $L_2(D)$. The inner product is defined as
\[ (u,v)=\int_{D} \Big ( uv + \frac{\partial u}{\partial x}\frac{\partial v}{\partial x}+ \frac{\partial u}{\partial t}\frac{\partial v}{\partial t} \Big )\, dx\,dt. \]

\subsection{Multiphase Stefan Optimal Control Problem}\label{E:1:1b}

Following the well-known reformulation of the IMSP (see \cite{LSU, Oleinik}), we consider the transformation
\begin{equation}\label{vt}
v(x,t)=F(u(x,t)):=\int\limits_{u_1}^{u(x,t)}k(y)\,dy.
\end{equation}
Then $v^j=\int\limits_{u^1}^{u^j}k(y)\,dy$, $~v^1=0<\cdots<v^J$, and our conditions become:
\begin{align}
\beta(v)v_t-v_{xx}=f(x,t), &\qquad (x,t)\in D, v(x,t)\neq v^j, \label{bvpde}\\
v|_{x=\xi_j(t)}=v^j, &\qquad 0<t\leq T, \label{vphase} \\
[v]|_{x=\xi_j(t)}=0, &\qquad 0<t\leq T,\quad  j=\overline{1,J},\label{vcont} \\
[v_x]|_{x=\xi_j(t)} = \gamma_j\frac{d\xi_j}{dt}, &\qquad 0<t\leq T,\quad  j=\overline{1,J}, \label{vjump} \\
v(x,0)=\Phi(x)=\int\limits_{u^1}^{\phi(x)}k(y)\,dy,&\qquad0\leq x\leq\ell, \label{vphi}\\
v_x|_{x=0}=g(t),&\qquad 0<t\leq T,  \label{vg}\\
v_x|_{x=\ell} = p(t),&\qquad 0<t\leq T,  \label{vp}\\
v(\ell,t)=\Gamma(t)=\int\limits_{u^1}^{\nu(t)}k(y)\,dy,&\qquad0<t\leq T,\label{vgamma}
\end{align}
where
\begin{equation}\label{beta}
\beta(v)=\frac{\alpha(F^{-1}(v))}{k(F^{-1}(v))},
\end{equation}
and $F^{-1}$ is an inverse function of $F$. The function $\beta(v)$ is of similar type as $\alpha$ and $k$, so that it is a positive $C^1$ function on each segment $(F(-\infty),v^1],$ $[v^j,v^{j+1}], j=1,\ldots,J-1$ and $[v^{J},F(+\infty))$ with 1st type discontinuity at $v=v^j, ~j=\overline{1,J}$. Now, we can invoke a function $b(v)$ such that $b'(v)=\beta(v)$.  Our partial differential equation becomes
\begin{equation}\label{bpde}
\frac{\partial b(v)}{\partial t}-v_{xx}=f(x,t), \qquad (x,t)\in D, v(x,t)\neq v^j.
\end{equation}
Moreover, we're free to choose the jump of $b$ at the values $v=v^j$. We choose them in such a way that $[b(v)]|_{v=v^j}=\gamma_j$ so that upon integration by parts of (\ref{bpde}) over $D$, the integrals over the phase transition boundaries cancel out. Define the level sets
\begin{equation}\label{surfaces}\nonumber
\n C_j:=\{(x,t)\in D~|~v(x,t)=v^j\},\qquad j=\overline{1,J}.
\end{equation} 

\textbf{Definition.} We say that a measurable function $B(x,t,v)$ is \textit{of type }$\n B$ if
\begin{enumerate}[(a)]
	\item $B(x,t,v) = b(v), \qquad v\neq v^j, \quad\forall j=\overline{1,J}$,$~$ and
	\item $B(x,t,v) \in [b(v^j)^-,b(v^j)^+],\qquad v=v^j$ for some $j$.
\end{enumerate}
Note that $B(x,t,v)$ can take different values for different $(x,t)$ when $v=v^j$ for some $j$.

Given $g$, a solution to the Stefan problem (\ref{bvpde})-(\ref{vp}) is understood in the following sense:\\

\textbf{Definition.} $v\in W_2^{1,1}(D)\cap L_{\infty}(D)$ is called a \emph{weak solution of the Stefan problem}  (\ref{bvpde})-(\ref{vp}) if for any two functions $B,B_0$ of type $\n B$, the following integral identity is satisfied:
\begin{gather}
\int\limits_0^T\int\limits_0^{\ell}\Big[-B(x,t,v(x,t))\psi_t+v_x\psi_x-f\psi\Big]\,dxdt - \int\limits_0^{\ell}B_0(x,0,\Phi(x))\psi(x,0)\,dx - \nonumber \\- \int\limits_0^Tp(t)\psi(\ell,t)\,dt+\int\limits_0^Tg(t)\psi(0,t)\,dt=0, \qquad \forall \psi\in W_2^{1,1}(D),\quad\psi(x,T)=0. \label{weaksol}
\end{gather}
Note that the definition of the weak solution doesn't impose any restriction on the structure of the free boundaries, which are identified as level sets of the weak solution. In particular, the level sets $\n C_j$ need not be curves in $D$, and the Stefan condition (\ref{stefancond}) is not explicitly required, but is absorbed into the definition of the weak solution through the function $B(x,t,v)$ and the integral identity (\ref{weaksol}).  \\

Consider the control set
\[\n G_R =\big\{g: g\in W_2^1(0,T),\,\,\Vert g\Vert_{W_2^1[0,T]}\leq R\big\}.\]

We wish to minimize the cost functional $\n J$ given by
\begin{equation}\label{functional}
\n J(g) = \Vert v(\ell,t;g)-\Gamma(t)\Vert^2_{L_2(0,T)}
\end{equation}
on $\n G_R$, where $v=v(x,t;g)\in W_2^{1,1}(D)\cap L_{\infty}(D)$ is a weak solution of the Stefan problem in the sense of (\ref{weaksol}). This optimal control problem will be called \emph{Problem $\m I$}.

\subsection{Discrete Optimal Control Problem}\label{E:1:1c}
Let 
\[\omega_{\tau} = \{t_k,k=\overline{1,n}\},\tau=\frac Tn, t_k=k\tau,\qquad\quad \omega_h= \{x_i,i=\overline{1,m}\},h=\frac{\ell}m, x_i=ih
\]
be grids in the time and space domains, respectively, and we'll assume from here on that 
\[
m\rightarrow\infty\quad\text{as } n\rightarrow\infty. \nonumber
\]
Define the Steklov averages
\begin{gather}\label{steklov}
a_k = \frac1{\tau}\int\limits_{t_{k-1}}^{t_k}a(t)\,dt,\qquad\qquad \Phi_i = \frac 1h\int\limits_{x_i}^{x_{i+1}}\Phi(x)\,dx,~~ \Phi_m=\Phi(\ell), \\ f_{ik}=\frac{1}{\tau h}\int\limits_{t_{k-1}}^{t_k} \int\limits_{x_i}^{x_{i+1}}f(x,t)\,dxdt,\qquad k=\overline{1,n},\quad i=\overline{0,m-1} \nonumber,
\end{gather}
where $a$ stands for any of the functions $p$, $\Gamma$, $g$, or $g^n$. Introduce the discretized control set
\[\n G_R^n = \{[g]_n\in\bb R^{n+1}: \Vert[g]_n\Vert_{w_2^1}\leq R\}\]
where $[g]_n=(g_0,g_1,\ldots,g_n)$, and
\[
\Vert[g]_n\Vert_{w_2^1}^2=\sum\limits_{k=1}^{n}\tau g_k^2+\sum\limits_{k=1}^n\tau g_{k\bar t}^2
\]
with $g_{k\bar t}=\frac{g_k-g_{k-1}}{\tau}.$ Assume that every element $g\in W_2^1(0,T)$ is continued to $[-\tau, 0]$ as a constant $g(0)$. Consider now the mappings between the discrete and continuous control sets, $\n Q_n: W_2^1(0,T)\rightarrow \bb R^{n+1},\quad \n P_n:\bb R^{n+1}\rightarrow W_2^1(0,T)$ as
\begin{gather}
\n Q_n(g)=[g]_n,\quad\text{for }g\in\n G_R,\,\text{where }g_k=\frac1{\tau}\int\limits_{t_{k-1}}^{t_k}g(t)\,dt,\quad k=\overline{0,n},\\
\n P_n([g]_n)=g^n,~\text{for }[g]_n\in\n G_R^n;~~ g^n(t)=g_{k-1}+\frac{g_k-g_{k-1}}{\tau}(t-t_{k-1}), t\in[t_{k-1},t_k),\,\, k=\overline{1,n}. \label{Pmap}
\end{gather}
Approximate the function $b(v)$ by the infinitely differentiable sequence 
\begin{equation}\label{smoothing}
b_{n}(v)=\int_{v-\frac{1}{n}}^{v+\frac{1}{n}}b(y)\omega_n(v-y)dy,
\end{equation}
where $\omega_n$ is a standard mollifier defined as
\begin{equation}\label{kernel}
\omega_n(v) =\left\{\begin{matrix}\m C n e^{-\frac{1} {1-n^2v^2}},\quad&|v|\leq\frac{1}{n}\\[2mm] 0,\quad&|v|>\frac{1}{n}\end{matrix}\right.
\end{equation}
and the constant $\m C$ is chosen so that $\int\limits_{\bb R}\omega_1(u)\,du =1$. Since $b'(v)$ is piecewise-continuous, we also have
\begin{equation}\label{smoothing_derivative}
b_{n}'(v)=\int_{v-\frac{1}{n}}^{v+\frac{1}{n}}b'(y)\omega_n(v-y)dy.
\end{equation}
Hence $b_n$ is strictly monotonically increasing. Next we define a discrete state vector, which represent the solution of the discrete multiphase Stefan problem.\\

\textbf{Discrete State Vector.} Given $[g]_n$, the vector function $[v([g]_n)]_n=\big(v(0),v(1),\ldots,v(n)\big);$ $\quad v(k)\in\bb R^{m+1},\quad k=0,\ldots,n$ is called a \emph{discrete state vector} if
\begin{enumerate}[(a)]
	\item $v_i(0) = \Phi_i,\quad i =\overline{0,m}, $\\
	\item For arbitrary $k=1,\ldots,n$, the vector $v(k)\in\bb R^{m+1}$ satisfies
	\begin{equation}\label{dsvsum}
	\sum\limits_{i=0}^{m-1}h\Big[\big(b_n(v_i(k))\big)_{\bar t} \eta_i+v_{ix}(k)\eta_{ix}-f_{ik}\eta_i\Big]-p_k\eta_m+g_k^n\eta_0=0,\qquad\forall \eta=(\eta_i)\in\bb R^{m+1}.
	\end{equation}
\end{enumerate}
Given $[g]_n\in\n G^n_R$, the discrete cost functional $\n I_n$ is defined as
\begin{equation}
\n I_n([g]_n) = \sum\limits_{k=1}^n\tau\Big(v_m(k)-\Gamma_k\Big)^2 \label{In}
\end{equation}
where $v_m(k)$ are components of the discrete state vector $[v([g]_n)]_n$. We define 
\[\n I_{n_*}:= \inf\limits_{[g]_n\in\n G_R^n}\n I_n([g]_n).\]
The discrete optimal control problem will be labeled \emph{Problem $\m I_n$}. Furthermore, the following interpolations will be considered:
\begin{gather}
\tilde v(x,t) = v_i(k),\qquad x\in[x_i,x_{i+1}],\quad t\in[t_{k-1},t_k],\qquad i =\overline{0,m-1},\quad k =\overline{0,n}, \nonumber \\ 
\hat v(x;k) = v_i(k)+v_{ix}(k)(x-x_i),\qquad x\in[x_i,x_{i+1}],\quad i=\overline{0,m-1}, \nonumber \\
v^{\tau}(x,t) = \hat v(x;k),\qquad t\in[t_{k-1},t_k], \nonumber \\
\hat v^{\tau}(x,t) = \hat v(x;k-1) +\hat v_{\bar t}(x;k)(t-t_{k-1}),\qquad t\in[t_{k-1},t_k],\quad k = \overline{1,n}. \label{interpolations}
\end{gather}

\subsection{Formulation of the Main Results}\label{E:1:1d}
Unless otherwise stated, throughout the paper we suppose that
\begin{gather}
f\in L_{\infty}(D),\quad p\in W_2^1(0,T),\quad\Phi\in W_2^1(0,\ell)
\end{gather}
and the one-dimensional Lebesgue measure of the level sets $\xi_j(0), j=\overline{1,J}$ defined in (\ref{initialinterfaceset}) is zero. Equivalently, the level sets $\{x\in[0,\ell]~|~ \Phi(x)=v^j\}$ are of Lebesgue measure null sets.
Concerning the behavior of the coefficients $\alpha$ and $k$ at $\infty$, we take the following assumptions:
\begin{gather}
\int_{u_1}^{\infty} k(u)du=\infty, \label{k_at_infinity}\\
\liminf_{u\to\infty} \frac{\alpha(u)}{k(u)} \ge a_0>0. \label{b_at_infinity}
\end{gather}
Condition \eqref{k_at_infinity} guarantees that the domain of $b(v)$ is $\bb R$. Condition \eqref{b_at_infinity} implies that there is a uniform positive lower bound
for $b'$ and $b_n'$:
\begin{equation}\label{bbar}
b'(v), b_n'(v) 
\geq \bar{b}>0, \ v\in \bb R.
\end{equation}
for some $\bar{b}>0$.
\begin{theorem}\label{existence}
	The Problem $\m I$ has a solution, that is, the set 
	\[\n G_*=\Big\{g\in\n G_R\big|\n J(g)=\n J_*:=\underset{g\in \n G_R}{\inf }\n J(g)\Big\}\]
	is not empty.
\end{theorem}
The proof of Theorem \ref{existence} hinges upon showing the weak continuity of the cost functional $\n J$. The weak continuity of  $\n J$ will be established by proving an $L_{\infty}(D)$ bound  and a $W_2^{1,1}(D)$ - energy estimation for the solution to the discrete Stefan problem, and subsequent use of compactness of the family of interpolations.

\begin{theorem}\label{convergence}
	The sequence of discrete optimal control problems $\m I_n$ approximates the optimal control problem $\m I$ with respect to functional, that is,
	\begin{equation}\label{Eq:W:1:18}
	\lim\limits_{n\to +\infty} \n{I}_{n_*}=\n{J}_*, 
	\end{equation}
	where
	\[ \n{I}_{n_*}=\inf\limits_{\n G_R^n} \n{I}_n([g]_n), \ n=1,2,\ldots ~.\]
	If $[g]_{n_\ep}\in \n G_R^n$ is chosen such that
	\[ \n{I}_{n_*} \le \n{I}_n([g]_{n_\ep})\le \n{I}_{n_*}+\ep_n, \ \ep_n \downarrow 0, \]
	then the sequence $g^n=\n{P}_n([g]_{n_\ep})$ has a subsequence convergent to some element $g_*\in\n G_*$ weakly
	in $W_2^1(0,T)$ and strongly in $L_2(0,T)$. Moreover, the piecewise linear interpolations $\hat{v}^\tau$ of the corresponding discrete state vectors $[v([g]_{n_\ep})]_n$ converge to the weak solution $v(x,t;g_*) \in W_2^{1,1}(D)$ of the Stefan Problem weakly in $W_2^{1,1}(D)$.
\end{theorem} 
The necessary and sufficient conditions for the convergence of discrete optimal control problems to the continuous optimal control problem are formulated in \cite{Vasilev} (see Lemma~\ref{Vasil} in Section~\ref{E:2:1} below). The proof of Theorem \ref{convergence} is based on the proof that the conditions of the general criteria are satisfied. As before, the $L_{\infty}$ bound  and the $W_2^{1,1}$ energy estimation for the solution to the discrete Stefan problem play a significant role in this context.  

\section{Preliminary Results}\label{E:2:1}
\begin{lemma}\label{eudsv} Given any $[g]_n\in\n G^n$, and any $h,\tau$, a discrete state vector exists uniquely.
\end{lemma}
\noindent\textit{Proof.} First we prove uniqueness. Suppose $v$ and $\tilde v$ both are discrete state vectors for a given $[g]_n$. Due to (a) from the discrete state vector definition, we have that $v(0)=\tilde v(0)$. For a fixed $k\geq 1$, suppose that $v(k-1)=\tilde v(k-1)$. Because (\ref{dsvsum}) is satisfied for both $v$ and $\tilde v$, subtract the identities for $\eta=v(k)-\tilde v(k)$ to get:
\[
\sum\limits_{i=0}^{m-1}\Big[\big(b_n(v_i(k))_{\bar t}-b_n(\tilde v_i(k))_{\bar t}\big)\big(v_i(k)-\tilde v_i(k)\big)+\big(v_{ix}(k) -\tilde v_{ix}(k)\big)^2\Big]=0.
\]
However,
\begin{align*}
b_n(v_i(k))_{\bar t}-b_n(\tilde v_i(k))_{\bar t} &= \frac{b_n(v_i(k))-b_n(v_i(k-1))}{\tau} - \frac{b_n(\tilde v_i(k))-b_n(\tilde v_i(k-1))}{\tau} \\[4mm] &= \frac{b_n(v_i(k))- b_n(\tilde v_i(k))}{\tau},
\end{align*}
so that the previous summation identity becomes:
\[
\sum\limits_{i=0}^{m-1}\Big[\frac1{\tau}\big(b_n(v_i(k))- b_n(\tilde v_i(k))\big)\big(v_i(k)-\tilde v_i(k)\big)+\big(v_{ix}(k) -\tilde v_{ix}(k)\big)^2\Big]=0.
\]
Since $b_n(v)$ is monotonically increasing, the whole summand is non-negative. Therefore, it is equal to 0, which  implies that $v_i(k)=\tilde v_i(k),\quad \forall i=\overline{0,m}$. Hence, by induction, $v=\tilde v$.

Now we seek to prove existence.  Again we'll rely on induction. Construct $v(0)$ as given in (a) of the Discrete State Vector Definition. Note that $\Vert v(0)\Vert:=\max\limits_i|v_i(0)| = \max\limits_i|\Phi_i|\leq \Vert\Phi\Vert_{L_{\infty}[0,\ell]}$.  Now fix $k\geq 1$, and assume that $v(k-1)$ has been constructed successfully so that (\ref{dsvsum}) is satisfied for all $K<k$. Moreover, assume that $\Vert v(k-1)\Vert <+\infty$. Notice that the summation identity (\ref{dsvsum}) is equivalent to solving the following system of non-linear equations:
\begin{equation}\label{system}
\left\{\begin{matrix}\left(v_0(k)+\frac{h^2}{\tau}b_n(v_0(k))\right)-v_1(k)&=&\frac{h^2}{\tau}b_n(v_0(k-1))+h^2f_{0k}-hg_k^n\\[4mm]
-v_{i-1}(k)+\left(2v_i(k)+\frac{h^2}{\tau}b_n(v_i(k))\right)-v_{i+1}(k) &=& \frac{h^2}{\tau}b_n(v_i(k-1))+h^2f_{ik},\quad i =\overline{1,m-1} \\[4mm] -v_{m-1}(k)+v_m(k)&=&hp_k \end{matrix}\right..
\end{equation}
We will construct $v(k)$ by the method of successive approximations. It is critical to remember that $h,\tau$ will be fixed here. Choose $v^0=v(k-1)$. Having obtained $v^N$,we search $v^{N+1}$ as a solution of the following system:
\begin{equation}\label{system1}
\left\{\begin{matrix}\left(v_0^{N+1}(k)+\frac{h^2}{\tau}b_n(v_0^{N+1}(k))\right)-v_1^N(k)&=&\frac{h^2}{\tau}b_n(v_0(k-1))+h^2f_{0k}-hg_k^n\\[4mm]
-v_{i-1}^N(k)+2v_i^{N+1}(k)+\frac{h^2}{\tau}b_n(v_i^{N+1}(k))-v_{i+1}^N(k) &=& \frac{h^2}{\tau}b_n(v_i(k-1))+h^2f_{ik},~~ i =\overline{1,m-1} \\[4mm] -v_{m-1}^{N+1}(k)+v_m^{N+1}(k)&=&hp_k \end{matrix}\right..
\end{equation}
We now proceed to prove that the sequence $\{v^N\}$ converges to the unique solution of (\ref{system}). Substract (\ref{system1}) for $N$ and $N-1$ to get
\begin{gather*}
\left\{\begin{matrix}v_0^{N+1}(k)-v_0^N(k)+\frac{h^2}{\tau}\left[b_n(v_0^{N+1}(k))- b_n(v_0^N(k))\right]=v_1^N(k)-v_1^{N-1}(k)\\[4mm]
2(v_i^{N+1}(k)-v_i^N(k))+\frac{h^2}{\tau}\left[b_n(v_i^{N+1}(k))-b_n(v_i^N(k)) \right]=v_{i+1}^N(k)-v_{i+1}^{N-1}(k)+v_{i-1}^N(k)-v_{i-1}^{N-1}(k)\\[4mm] v_m^{N+1}(k)-v_m^N(k)=v_{m-1}^{N+1}(k)-v_{m-1}^N(k) \end{matrix}\right.
\end{gather*}
which is transformed as
\begin{gather*}
\left\{\begin{matrix}\left(1+\frac{h^2}{\tau}\zeta_{n,N}^0\right)(v_0^{N+1}(k)- v_0^N(k))=v_1^N(k)-v_1^{N-1}(k)\\[4mm]
\left(2+\frac{h^2}{\tau}\zeta_{n,N}^i\right)(v_i^{N+1}(k)-v_i^N(k))=(v_{i+1}^N(k)- v_{i+1}^{N-1}(k))+(v_{i-1}^N(k)-v_{i-1}^{N-1}(k))\\[4mm] v_m^{N+1}(k)-v_m^N(k)=v_{m-1}^{N+1}(k)-v_{m-1}^N(k) \end{matrix}\right.
\end{gather*}
where
\begin{equation*}
\zeta_{n,N}^i:=\int_0^1b_n'(\theta v_i^{N+1}(k)+(1-\theta)v_i^N(k))d\theta, \quad i =\overline{0,m-1}.
\end{equation*}
Due to \eqref{bbar}, we have $\zeta_{n,N}^i\geq\bar b,~~i=\overline{0,m-1}$. Hence we have
\begin{gather}\label{system2}
\left\{\begin{matrix}v_0^{N+1}(k)-v_0^N(k)=\frac{v_1^N(k)-v_1^{N-1}(k)}{ 1+\frac{h^2}{\tau}\zeta_{n,N}^0}\\[6mm]
v_i^{N+1}(k)-v_i^N(k)=\frac{v_{i+1}^N(k)- v_{i+1}^{N-1}(k)+v_{i-1}^N(k)-v_{i-1}^{N-1}(k)}{2+\frac{h^2}{\tau}\zeta_{n,N}^i}, \qquad i=\overline{1,m-1}\\[6mm] v_m^{N+1}(k)-v_m^N(k)=v_{m-1}^{N+1}(k)-v_{m-1}^N(k) \end{matrix}\right..
\end{gather}
Let $A_N:=\max\limits_{0\leq i\leq m} {\big|v_i^{N+1}(k)-v_i^N(k)\big|},\qquad\delta=\left(1+\frac{h^2}{2\tau}\frac{\bar b}{2}\right)^{-1}\in(0,1)$. From (\ref{system2}) it easily follows
\begin{equation}\label{A}
A_N\leq\delta A_{N-1}\leq\cdots\leq A_0\delta^N.
\end{equation}
Now it is possible to prove that there exist finite limits
\begin{equation}\label{limit}
v_i(k)=\lim\limits_{N\rightarrow+\infty}v_i^N(k),\qquad i=0,1,\ldots,m.
\end{equation}
Indeed, from (\ref{A}) it follows that for arbitrary $i=0,1,\ldots,m$, we have
\begin{equation}\label{bounds}
-A_0\delta^N\leq v_i^{N+1}(k)-v_i^N(k)\leq A_0\delta^N.
\end{equation}
By summation we have
\begin{equation}\label{sumbounds}
v_i^N(k)-A_0\sum\limits_{\ell=N}^{+\infty}\delta^{\ell}\leq v_i^N(K)-A_0 \sum\limits_{\ell=N}^{M-1}\delta^{\ell}\leq v_i^M(k)\leq v_i^N(k)+A_0 \sum\limits_{\ell=N}^{M-1}\delta^{\ell}\leq v_i^N+A_0 \sum\limits_{\ell=N}^{+\infty}\delta^{\ell}
\end{equation}
for all $M>N\geq0$. In particular, by choosing $N=0$ it follows that the sequence $\{v^N\}$ is bounded in $\bb R^{m+1}$. Let us now assume
\[
\liminf\limits_{N\rightarrow+\infty}v_i^N(k)=\lim\limits_{p\rightarrow+\infty} v_i^{N_p}(k),\qquad N_p<N_{p+1},~~p=0,1,\ldots;\quad\lim\limits_{p\rightarrow +\infty}N_p=+\infty.
\]
By choosing in (\ref{sumbounds}) $N=N_p$ we have
\[
v_i^M(k)\leq v_i^{N_p}(k)+A_0\sum\limits_{\ell=N_p}^{+\infty}\delta^{\ell},\qquad M>N_p
\]
which implies that
\[
\limsup\limits_{M\rightarrow\infty}v_i^M(k)\leq v_i^{N_p}(k)+A_0\sum\limits_{\ell=N_p }^{+\infty}\delta^{\ell},\qquad p=1,2,\ldots
\]
Passing to limit as $p\rightarrow+\infty$ we have
\[
\limsup\limits_{M\rightarrow\infty}v_i^M(k)\leq\lim\limits_{p\rightarrow+\infty} v_i^{N_p}(k)=\liminf\limits_{N\rightarrow+\infty}v_i^N(k).
\]
Since opposite inequality is obvious, it follows that finite limits (\ref{limit}) exist.\hfill{$\square$}

Given the existence and uniqueness of the discrete state vector for fixed $n$, we can uniquely define for each $k=1,\ldots,n$ the vector $\zeta_k$ whose $m$ components $\zeta_k^i$ are given by
\begin{equation}
\zeta^i_k=\int_0^1b_n'(\theta v_i(k)+(1-\theta)v_i(k-1))d\theta, \quad i =\overline{0,m-1}\label{zetav}.
\end{equation}

\begin{lemma}\label{Vasil}\cite{Vasilev} The sequence of discrete optimal control problems $\m I_{n}$ approximates the continuous optimal control problem $I$ if and only if the following conditions are satisfied:
	\begin{description}
		\item{\bf(1)} for arbitrary sufficiently small $\ep>0$ there exists $M_1=M_1(\ep)$ such that $\n{Q}_M(g)\in \n G^{M}_R$ for all $g \in \n G_{R-\ep}$ and $M\ge M_1$; and for any fixed $\ep>0$ and for all $g\in \n G_{R-\ep}$ the following inequality is satisfied:
		\begin{equation}\label{firstcondition}
		\limsup\limits_{M\to \infty} \Big ( \n{I}_M(\n{Q}_M(g))-\n{J}(g) \Big ) \le 0.
		\end{equation}
		\item{\bf(2)} for arbitrary sufficiently small $\ep>0$ there exists  $M_2=M_2(\ep)$ such that $\n{P}_M([g]_{M})\in \n G_{R+\ep}$ for all $[g]_{M} \in \n G^{M}_R$ and $M\ge M_2$; and for all $[g]_{M}\in \n G^{M}_R$, $M\ge 1$ the following inequality is satisfied:
		\begin{equation}\label{secondcondition}
		\limsup\limits_{M\to \infty} \Big ( \n{J}(\n{P}_M([g]_{M})) -\n{I}_M([g]_{M})  \Big ) \le 0.
		\end{equation}
		\item{\bf(3)} the following inequalities are satisfied:
		\begin{align}
		\limsup\limits_{\ep \to 0} \n{J}_*(\ep) \ge \n{J}_*, 
		\ \ \liminf\limits_{\ep \to 0} \n{J}_*(-\ep) \le \n{J}_*,
		\end{align}\label{thirdcondition}
		where $\n{J}_*(\pm\ep)=\inf\limits_{\n G_{R\pm \ep}}\n{J}(g)$.
	\end{description}
\end{lemma}

\begin{lemma}\label{mappings} The mappings $\n P_n, \n Q_n$ satisfy the conditions of Lemma \ref{Vasil}.
\end{lemma}
\noindent\emph{Proof. } Let $g\in\n G_{R}$, $[g]_n=\n Q_n(g)$. We observe that
\begin{gather}
\sum\limits_{k=1}^{n}\tau g_k^2 =\sum\limits_{k=1}^n\tau\left(\frac1{\tau}\int\limits_{t_{k-1}}^{t_k}g(t)\,dt \right)^2 \leq \int\limits_0^Tg^2(t)\,dt = \Vert g\Vert_{L_2[0,T]}^2, \label{q1}\\[5mm]
\sum\limits_{k=1}^n\tau g_{k\bar t}^2=\sum\limits_{k=1}^n\frac1{\tau}\left(\frac1{\tau}\int\limits_{t_{k-1}}^{t_k} g(t)\,dt-\frac1{\tau}\int\limits_{t_{k-2}}^{t_{k-1}} g(t)\,dt\right)^2 \leq \sum\limits_{k=1}^n\frac1{\tau^3}\left(\int\limits_{t_{k-1}}^{t_k} \int\limits_{t-\tau}^tg'(\xi)\,d\xi dt\right)^2 \nonumber \\\leq \sum\limits_{k=1}^n\frac1{\tau}\int\limits_{t_{k-1}}^{t_k} \int\limits_{t-\tau}^t|g'(\xi)|^2\,d\xi dt
\leq \Vert g'\Vert_{L_2[0,T]}^2. \label{q2}
\end{gather}
From (\ref{q1}),(\ref{q2}), we get
\begin{equation}\label{q3}
\Vert[g]_n\Vert_{w_2^1}^2 \leq R^2.
\end{equation}
Now let $[g]_n\in\n G_R^n$ be given and write $g^n=\n P_n([g]_n)$. We see that
\begin{gather}
\int\limits_0^T\left|\frac{dg^n(t)}{dt}\right|^2\,dt = \sum\limits_{k=1}^n\,\int\limits_{t_{k-1}}^{t_k}g_{k\bar t}^2\,dt = \sum\limits_{k=1}^n\tau g_{k\bar t}^2 ,\label{p1} \\[5mm]
\int\limits_0^T|g^n(t)|^2\,dt =\sum\limits_{k=1}^n\, \int\limits_{t_{k-1}}^{t_k} \Big(g_{k-1}+g_{k\bar t}(t-t_{k-1})\Big)^2\,dt\nonumber \\ = \sum\limits_{k=1}^n\tau g_{k-1}^2 + \sum\limits_{k=1}^n\tau^2g_{k-1}g_{k\bar t}+ \frac13\sum\limits_{k=1}^n\tau^3g_{k\bar t}^2 =\sum\limits_{k=1}^n\tau g_kg_{k-1} + \frac13\sum\limits_{k=1}^n\tau^3g_{k\bar t}^2 .\label{p2}
\end{gather}
Actually, since $[g]_n\in\n G_R^n$, it is the case that $\sum\limits_{k=1}^n\tau g_{k\bar t}^2 \leq C^2$ where $C$ is a constant independent of $n$. This implies that
\begin{equation}\label{gnlipschitz}
|g_k-g_{k-1}|\leq C\sqrt{\tau},\qquad k=1,\ldots,n.
\end{equation}
Using (\ref{gnlipschitz}) in (\ref{p2}), we can write
\begin{gather}
\int\limits_0^T|g^n(t)|^2\,dt\leq\sum\limits_{k=1}^n\tau g_k^2 +C\sqrt{\tau}\sum\limits_{k=1}^n\tau|g_k| + \frac13\tau^2\sum\limits_{k=1}^n\tau g_{k\bar t}^2 \nonumber \\ \leq \sum\limits_{k=1}^n\tau g_k^2 + C\sqrt T \sqrt{\tau}\sqrt{\sum\limits_{k=1}^n\tau g_k^2} + \frac13\tau^2\sum\limits_{k=1}^n\tau g_{k\bar t}^2 .\label{p3}
\end{gather}
Combining (\ref{p1}) and (\ref{p3}) we have that
\begin{equation}
\Vert g^n\Vert_{W_2^1[0,T]}^2 \leq \Vert[g]_n\Vert_{w_2^1}^2 +C\sqrt T \sqrt{\tau}\sqrt{\sum\limits_{k=1}^n\tau g_k^2} + \frac13\tau^2\sum\limits_{k=1}^n\tau g_{k\bar t}^2 \leq R^2+O(\sqrt{\tau})  \label{p4}.
\end{equation}
Owing to (\ref{p4}), we can choose $n$ so large that $\tau$ will be small enough to guarantee that the right-hand side will be bounded by $(R+\ep)^2$. Hence $g^n\in\n G_{R+\ep}$ for all $n$ large enough.\hfill{$\square$}

\begin{lemma}\label{unique} There is at most one solution to the Stefan problem in the sense of (\ref{weaksol}). 
\end{lemma}

\noindent\emph{Proof.} That a solution to the Stefan problem in the sense of (\ref{weaksol}) is unique follows by an argument analogous to that presented in Section $9$ of Chapter V of \cite{LSU}. Indeed, we will prove uniqueness in a wider class of solutions than that given in (\ref{weaksol}). Suppose that $v\in L_{\infty}(D)$ only, not necessarily in the Sobolev space $W_2^{1,1}(D)$, and that for any two functions $B,B_0$ of type $\n B$ it satisfies the identity
\begin{gather}
\int\limits_0^T\int\limits_0^{\ell}\Big[B(x,t,v)\psi_t+v\psi_{xx}+f\psi\Big]\,dxdt + \int\limits_0^{\ell}B_0(x,0,\Phi(x))\psi(x,0)\,dx + \nonumber \\ \int\limits_0^Tp(t)\psi(\ell,t)\,dt-\int\limits_0^Tg(t)\psi(0,t)\,dt=0, \qquad \forall \psi\in W_2^{2,1}(D), \psi(x,T)=0, \psi_x(0,t)=\psi_x(\ell,t)=0. \label{weakersol}
\end{gather}
The class of functions satisfying the above definition contains the class of solutions given in (\ref{weaksol}). Suppose $v$ and $\tilde v$ are two solutions in the sense of (\ref{weakersol}). Due to our assumption on $\Phi$, subtracting (\ref{weakersol}) with $\tilde v$ from that with $v$ guarantees that the second integral in (\ref{weakersol}) vanishes, and we obtain:
\[
\int\limits_0^T\int\limits_0^{\ell}\big(B(x,t,v)-\tilde B(x,t,\tilde v)\big)\left(\psi_t+a(x,t)\psi_{xx} \right)\,dx\,dt = 0
\]
where $a(x,t) = \frac{v-\tilde v}{B(x,t,v)-\tilde B(x,t,\tilde v)}$. For $(x,t)\in D$ such that $v(x,t)=\tilde v(x,t)$, it is the case that $a(x,t)=0$. Otherwise, since $B$ and $\tilde B$ are strictly increasing on $v$ a.e. $(x,t)\in D$, it follows that $a$ is non-negative for a.e. $(x,t)$. Moreover, the a.e. positiveness of $b'(v(x,t))$ implies that $\bar b =\text{essinf}~ b' > 0$ and that $b$ is strictly increasing, and so for almost every $(x,t)$ (assume that $\tilde v(x,t) < v(x,t)$ for the sake of notational simplicity),
\begin{align*}
|a(x,t)| & = \left|\frac{v-\tilde v}{\int\limits_{\tilde v(x,t)}^{v(x,t)}b'(w)\,dw+\sum\limits_{i: v^i\in(\tilde v(x,t),v(x,t))}(b(v^{i})^+-b(v^{i})^-)}\right| \leq\left|\frac{v-\tilde v}{\int\limits_{\tilde v}^v\bar b\,dv}\right| = \frac1{\bar b},
\end{align*}
so that $a$ is essentially bounded, and $\text{esssup}~ a(x,t) = a_1 < +\infty$. Fix $\ep>0$, and take as $\psi(x,t)$ the solution of the Neumann problem
\begin{equation}
\psi_t+(a(x,t)+\ep)\psi_{xx}=F(x,t),\qquad \psi_x(0,t)=\psi_x(\ell,t)= 0,\quad \psi(x,T)=0 \label{conjugate},
\end{equation}
where the $\ep$ is added to ensure the uniform parabolicity of the conjugate diffusion coefficient, and $F$ is an arbitrary smooth bounded function in $D$. Note that (\ref{conjugate}) is the conjugate heat equation. There exists a unique solution $\psi^{\ep} \in W_2^{2,1}(D)$ of the problem (\ref{conjugate})(\cite{LSU}). Our goal here is to use the arbitrariness of $F$ to obtain that $B-\tilde B=0$ a.e.; to this end, notice that through the use of (\ref{conjugate}), we can write
\begin{equation}
\int\limits_0^T\int\limits_0^{\ell}\big(B(x,t,v)-\tilde B(x,t,\tilde v)\big)\left(F-\ep\psi^{\ep}_{xx} \right)\,dx\,dt = 0 \label{almost}.
\end{equation}

Thus our goal will be attained if we have an energy estimate on $\psi_{xx}$ for solutions of (\ref{conjugate}). In the following, we prove a sufficient estimation for the analogous Heat Equation (the result follows immediately for the conjugate one by a simple change of variables). Let $a^{\ep}(x,t) =a(x,t)+\ep$, and for simplicity we don't write the superscript. Multiply the non-conjugate version of (\ref{conjugate}) by $\psi_{xx}$ and integrate it over the rectangle $D_t :=(0,\ell)\times(0,t)$ to get
\begin{gather}
-\int\limits_0^t\int\limits_0^{\ell}(\psi_{\tau}-a\psi_{xx})\psi_{xx}\,dx\,d\tau = -\int\limits_0^t\int\limits_0^{\ell}F\psi_{xx} \,dx\,d\tau= \int\limits_0^t\int\limits_0^{\ell}F_x\psi_x\,dx\,d\tau-\int\limits_0^tF\psi_x\Big|_0^{\ell}\,d\tau, \nonumber \\
\int\limits_0^t\int\limits_0^{\ell}\big((\psi_{\tau})_x\psi_x+a\psi^2_{xx}\big)\,dx\,d\tau- \int\limits_0^t\psi_{\tau}\psi_x \Big|_0^{\ell} \,d\tau = \int\limits_0^t\int\limits_0^{\ell}\big((\psi_x)_{\tau}\psi_x+a\psi^2_{xx}\big)\,dx\,d\tau = \int\limits_0^t\int\limits_0^{\ell}F_x\psi_x\,dx\,d\tau ,\nonumber \\
\frac12\int\limits_0^{\ell}\psi_x^2(x,t)\,dx-\frac12\int\limits_0^{\ell}\psi_x^2(x,0)\,dx +\int\limits_0^t\int\limits_0^{\ell}a \psi_{xx} ^2 \,dx\,d\tau \leq \frac12\int\limits_0^t\int\limits_0^{\ell}(\psi_x^2+F_x^2)\,dx\,d\tau, \nonumber \\
\int\limits_0^{\ell}\psi_x^2(x,t)\,dx +2\int\limits_0^t\int\limits_0^{\ell}a\psi_{xx}^2\,dx\,d\tau \leq \int\limits_0^t\int\limits_0^{\ell}\psi_x^2\,dx\,d\tau +\int\limits_0^t\int\limits_0^{\ell}F_x^2\,dx\,d\tau .\label{ue1}
\end{gather}
Letting now $y(t) = \int\limits_0^t\int\limits_0^{\ell}\psi_x^2\,dx\,d\tau$, it is clear that $y'(t) = \int\limits_0^{\ell}\psi_x^2(x,t) \,dx$, thus that (\ref{ue1}) implies
\[
y'(t)\leq y(t)+\int\limits_0^t\int\limits_0^{\ell}F_x^2\,dx\,d\tau.
\]
By Gronwall's Inequality now (more precisely Lemma 5.5 from \cite{LSU}), we deduce from the above differential inequality that $y(t) \leq \Big[e^t-1\Big]\int\limits_0^t\int\limits_0^{\ell}F_x^2\,dx\,d\tau$, or in other words
\begin{gather*}
\int\limits_0^t\int\limits_0^{\ell}\psi_x^2(x,\tau) \,dx\,d\tau \leq \Big[e^t-1\Big]\int\limits_0^t\int\limits_0^{\ell}F_x^2\,dx\,d\tau, \qquad \forall t\in(0,T],\intertext{so that by (\ref{ue1}),}\int\limits_0^{\ell}\psi_x^2(x,t) \,dx + 2\int\limits_0^t\int\limits_0^{\ell}a\psi_{xx}^2 \,dx\,d\tau \leq e^t\int\limits_0^t\int\limits_0^{\ell}F_x^2\,dx\,d\tau,\qquad \forall t\in(0,T].
\end{gather*}
The first of the above inequalities implies that $\underset{0\leq t\leq T}{\text{ess sup}}\int\limits_0^{\ell}\psi_x^2(x,t) \,dx\leq e^T\int\limits_0^T\int\limits_0^{\ell}F_x^2\,dx\,d\tau$. Now, since $\psi_t = a\psi_{xx}+F$, we have
\begin{align*}
\Vert\psi_t\Vert^2_{L_2(D_t)} &= \Vert a\psi_{xx}+F\Vert^2_{L_2(D_t)} \leq (\Vert a\psi_{xx}\Vert_{L_2(D_t)}+\Vert F\Vert_ {L_2(D_t)})^2 \leq 2\Vert a\psi_{xx}\Vert^2_{L_2(D_t)}+2\Vert F\Vert^2_{L_2(D_t)} \\ 
&\leq2\Vert F\Vert^2_{L_2(D_t)} +2a_0\int\limits_0^t\int\limits_0^{\ell}a\psi_{xx}^2\,dx\,d\tau \leq 2\Big(\Vert F\Vert^2_{L_2(D_t)}+a_0e^t\Vert F_x\Vert^2_{L_2(D_t)}\Big).
\end{align*}
These results combined provide the energy estimate we need:
\begin{equation}
\int\limits_0^T\int\limits_0^{\ell}\Big(\psi_t^2+a(x,t)\psi_{xx}^2\Big)\,dx\,dt+\underset{0\leq t\leq T}{\text{ess sup}}~\Vert\psi_x \Vert_{L_2(0,\ell)}^2 \leq 2\Big(\Vert F\Vert^2_{L_2(D)}+a_0e^T\Vert F_x\Vert^2_{L_2(D)}\Big). \label{eenp}
\end{equation}
Having (\ref{eenp}), we can now observe that
\begin{align*}
\left|\int\limits_0^T\int\limits_0^{\ell}(B-\tilde B)\ep\psi^{\ep}_{xx}\,dx\,dt\right| &= \left|\int\limits_0^T\int\limits_0^{\ell}(B-\tilde B)\frac{\ep}{(a+\ep)^{\frac12}}(a+\ep)^{\frac12}\psi^{\ep}_{xx}\,dx\,dt\right| \\
& \leq 2\text{esssup}~b(v) \left(\int\limits_0^T\int\limits_0^{\ell}\frac{\ep^2}{(a+\ep)}\,dx\,dt\right)^{\frac12} \left(\int\limits_0^T\int\limits_0^{\ell}(a+\ep)(\psi^{\ep}_{xx})^2\,dx\,dt\right)^{\frac12}\\
&\leq C\sqrt{\ep}\text{esssup}~b(v) \left(\int\limits_0^T\int\limits_0^{\ell}\frac{\ep}{(a+\ep)}\,dx\,dt\right)^{\frac12} \Vert F\Vert_{W_2^{1,0}(D)} \rightarrow 0 ~~\text{as}~~\ep\rightarrow0
\end{align*}
where $C$ is a constant depending only on $T$ and $a_0$. Recall that $\ep\leq a+\ep$, and so the integral on the right-hand side of the above inequality is bounded above by the area of the rectangle $D$. Therefore, (\ref{almost}) now implies
\[
\int\limits_0^T\int\limits_0^{\ell}\big(B(x,t,v(x,t))-\tilde B(x,t,\tilde v(x,t))\big)F\,dx\,dt = 0.
\]
Owing to the arbitrariness of $F$, the above equality gives that $B(x,t,v(x,t)) = \tilde B(x,t,\tilde v(x,t))$ a.e.$(x,t)\in D$, implying $b(v(x,t))=b(\tilde v(x,t)), ~\text{a.e. }(x,t)~$s.t.$~v(x,t)\neq v^j, j=1,\ldots,m$. Since $b$ is strictly increasing, we therefore have $v(x,t)=\tilde v(x,t)\quad\text{a.e.}~(x,t)$, so $v$ and $\tilde v$ coincide as solutions in the sense of (\ref{weakersol}), and thus we have proven uniqueness in this large class of solutions. \hfill{$\square$}

\begin{corollary}\label{Cjmeasure0} If a weak solution exists, all of the sets $\n C_j, j=\overline{1,J}$ have 2-dimensional measure 0.
\end{corollary}
\noindent\emph{Proof. } The proof of uniqueness gives us that $B_1(x,t,v(x,t))=B_2(x,t,v(x,t))$ a.e. on $D$, for \emph{any} two functions $B_1, B_2$ of type $\n B$. The functions of type $\n B$ generally differ on the sets $\n C_j$, so if one of them has positive measure, we arrive at a contradiction to Lemma \ref{unique}.\hfill{$\square$}

\section{Proofs of the Main Results}\label{E:3:1}

\subsection{$L_{\infty}$-estimation for the Discrete Stefan Problem}\label{E:3:1a}

\begin{theorem}\label{uniformboundedness} Suppose that $p\in L_{\infty}(0,T), \Phi\in L_{\infty}(0,\ell), f\in L_{\infty}(D)$. For $[g]_n\in \n G_R^n$ and $n,m$ large enough, the discrete state vector $[v([g]_n)]_n$ satisfies the following estimate:
	\begin{equation}\label{linfestimate}
	\Vert[v]_n\Vert_{\ell_{\infty}}:=\max\limits_{0\leq k\leq n}\Big(\max\limits_{0\leq i\leq m}|v_i(k)|\Big) \leq C_{\infty}\Big(\Vert f\Vert_{L_{\infty}(D)}+\Vert p\Vert_{L_{\infty}(0,T)}+\Vert g^n\Vert_{W_2^1(0,T)}+\Vert\Phi\Vert_{L_{\infty}(0,\ell)}\Big)
	\end{equation}
	where $C_{\infty}$ is a constant independent of $n$ and $m$.
\end{theorem}
\noindent\emph{Proof. } Fix $n$ arbitrarily large. Note $\max|v_i(0)|\leq\Vert\Phi\Vert_{L_{\infty}(0,\ell)}$. Consider a positive function $\gamma(x)\in C^2[0,\ell]$ satisfying
\begin{equation}\label{gamma}
\gamma(0)=\frac12,\qquad \gamma(\ell)=\frac12,\qquad \gamma'(0)=1,\quad\gamma'(\ell)=-1,\qquad \frac14\leq\gamma(x)\leq1,~x\in[0,\ell].
\end{equation}
Define $\gamma_i=\gamma(x_i),~i=\overline{0,m}$, and denote as $x^i$ the value in $[x_i,x_{i+1}]$ that satisfies (by MVT) $\gamma(x_{i+1})-\gamma(x_i)=\gamma'(x^i)h$. Transform the discrete state vector as
\[
w_i(k)=v_i(k)\gamma_i,\qquad i =\overline{0,m},\qquad k =\overline{0,n}.
\]
System (\ref{system}) can be rewritten as:
\begin{equation}\label{dpde}
\left\{\begin{matrix}h\zeta_k^0v_{0\bar t}(k)-v_{0x}(k)&=& hf_{0k}-g_k^n\\[4mm]
\zeta_k^iv_{i\bar t}(k)-v_{ix\bar x}(k) &=& f_{ik},\qquad i =1,\ldots, m-1 \\[4mm] v_{m-1,x}(k)&=&p_k \end{matrix}\right..
\end{equation}
We note
\begin{gather*}
v_i(k) = \frac1{\gamma_i}w_i(k),\qquad v_{i\bar t}(k) = \frac1{\gamma_i}w_{i\bar t}(k), \\[4mm]
v_{ix}(k) = \frac1{\gamma_{i+1}}w_{ix}(k)+\left(\frac1{\gamma_i}\right)_xw_i(k) = \frac1{\gamma_{i}}w_{ix}(k)+\left(\frac1{\gamma_i}\right)_xw_{i+1}(k), \\[4mm]
v_{ix\bar x}(k) = \frac1{\gamma_{i-1}}w_{ix\bar x} (k) +\left[ \left(\frac1{\gamma_i}\right)_{\bar x} + \left(\frac1{\gamma_i}\right)_{x} \right]w_{ix}(k)+ \left(\frac1{\gamma_i}\right)_{x\bar x}w_i(k) \\[2mm]\qquad\qquad\qquad =\frac1{\gamma_{i+1}}w_{ix\bar x} (k) +\left[ \left(\frac1{\gamma_i}\right)_{\bar x} + \left(\frac1{\gamma_i}\right)_{x} \right]w_{i\bar x}(k)+ \left(\frac1{\gamma_i}\right)_{x\bar x}w_i(k), \\[4mm]
\left(\frac1{\gamma_i}\right)_x = -\frac1{\gamma_i\gamma_{i+1}}\gamma_{ix},\qquad \left(\frac1{\gamma_i}\right)_{x\bar x} = -\frac1{\gamma_i\gamma_{i+1}}\gamma_{ix\bar x}+\frac{\gamma_{ix}+\gamma_{i\bar x}}{\gamma_{i-1}\gamma_i\gamma_{i+1}}\gamma_{i\bar x}.
\end{gather*}
Thus $w_i(0)=\gamma_i\Phi_i,\quad i=\overline{0,m}$, and for $k=\overline{1,n}$,
\begin{equation}\label{systemw}
\begin{matrix}\frac h{\gamma_0}\zeta_k^0w_{0\bar t}(k)-\frac1{\gamma_1}w_{0x}(k)-\left(\frac1{\gamma_0}\right)_xw_0(k)= hf_{0k}-g_k^n\\[4mm]
\frac1{\gamma_i}\zeta_k^iw_{i\bar t}(k)-\frac1{\gamma_{i-1}}w_{ix\bar x}(k) -\left[ \left(\frac1{\gamma_i}\right)_{\bar x} + \left(\frac1{\gamma_i}\right)_{x} \right]w_{ix}(k) - \left(\frac1{\gamma_i}\right)_{x\bar x}w_i(k)=f_{ik},\qquad i =\overline{1,m-1} \\[4mm] \frac1{\gamma_{m-1}}w_{m-1,x}(k)+\left(\frac{1}{\gamma_{m-1}}\right)_xw_m(k)=p_k \end{matrix}.
\end{equation}
Furthermore, transform $w_i(k)$ as:
\begin{equation}
\label{u} u_i(k) = w_i(k)e^{-\lambda t_k},\qquad i = \overline{0,m},\quad k =\overline{0,n}
\end{equation}
where
\begin{equation}\label{lambda}
\lambda = \frac{65}{\bar b} \big(\Vert\gamma''\Vert_{C[0,\ell]} +\Vert\gamma'\Vert_{C[0,\ell]}^2\big),
\end{equation}
and if $t^k\in[t_{k-1},t_k]$ satisfies through the MVT that $e^{\lambda t_k}-e^{\lambda t_{k-1}}=\lambda e^{\lambda t^k}(t_k-t_{k-1}) = \lambda e^{\lambda t^k}\tau$, then 
\[
w_{i\bar t}(k) = e^{\lambda t_{k-1}}u_{i\bar t}(k)+\lambda e^{\lambda t^k}u_i(k).
\] 
So $u_i(0) = w_i(0)=\gamma_i\Phi_i,\quad i=\overline{0,m}$, and for $k=\overline{1,n}$, the vector $u(k)$ satisfies the system
\begin{equation}
\label{systemu}\begin{matrix}\frac h{\gamma_0}\zeta_k^0e^{-\lambda\tau}u_{0\bar t}(k)-\frac1{\gamma_1}u_{0x}(k)+\left[\frac{\lambda h}{\gamma_0}\zeta_k^0e^{-\lambda(t_k-t^k)} -\left(\frac1{\gamma_0}\right)_x\right]u_0(k)= e^{-\lambda t_k}\big(hf_{0k}-g_k^n\big)\\[6mm]
\frac1{\gamma_i}\zeta_k^ie^{-\lambda\tau}u_{i\bar t}(k)-\frac1{\gamma_{i-1}}u_{ix\bar x}(k) -\left[ \left(\frac1{\gamma_i}\right)_{\bar x} + \left(\frac1{\gamma_i}\right)_{x} \right]u_{ix}(k) ~\qquad\qquad\qquad\qquad \\[3mm] \qquad\qquad\qquad +~ \left[\frac{\lambda}{\gamma_i}\zeta_k^ie^{-\lambda(t_k-t^k)} - \left(\frac1{\gamma_i}\right)_{x\bar x}\right]u_i(k)=f_{ik}e^{-\lambda t_k},\qquad i =\overline{1,m-1} \\[6mm] \frac1{\gamma_{m-1}}u_{m-1,x}(k)+\left(\frac{1}{\gamma_{m-1}}\right)_xu_m(k)=e^{-\lambda t_k}p_k \end{matrix}.
\end{equation} 

Now fix $k_1\leq n$, and define the following sets of indexes for convenience:
\begin{gather*}
\n M_{k_1} =\{(i,k)|i=0,\ldots,m,\quad k = 0,\ldots, k_1\} , \\
\n N =\{(i,k)|i=1,\ldots,m-1,\quad k = 1,\ldots, k_1\}, \\
\n T_0 = \{(i,k)|i=0,k = 1,\ldots, k_1\}, \\
\n T_m = \{(i,k)|i=m,k = 1,\ldots, k_1\}, \\
\n X_0 = \{(i,k)|i=0,\ldots,m,\quad k =0\}.
\end{gather*}
Unless confusion may arise, we omit the subscript to $\n M_{k_1}$. It is clear that 
\[
\n M = \n N\cup\n T_0\cup\n T_m\cup\n X_0.
\]
If $u_i(k)\leq 0$ in $\n M$, then $\max\limits_{\n M}u_i(k)\leq 0$. Suppose that there exists $(i,k)$ such that $u_i(k)>0$. Then $\max\limits_{\n M}u_i(k)>0$. Let $(i^*,k^*)\in \n M$ be such that $u_{i^*}(k^*) = \max\limits_{\n M}u_i(k)$.

If $(i^*,k^*)\in\n X_0$, then $u_{i^*}(k^*) = \max\limits_i\gamma_i\Phi_i \leq \max\limits_i\Phi_i \leq \max\limits_{[0,\ell]}\Phi(x)$.

If $(i^*,k^*)\in\n T_m$, then $i^*=m,~u_{m-1,x}(k^*)\geq0$ and we can choose $h$ small enough that $\gamma_{m-1,x}=\gamma'(x^{m-1})\in(-\frac32,-\frac12)$ so that
\[
-\frac{\gamma_{m-1,x}}{\gamma_m\gamma_{m-1}}u_m(k^*)\leq e^{-\lambda t_{k^*}}p_{k^*} \quad\implies\quad u_m(k^*)\leq \frac{\gamma_m\gamma_{m-1}}{-\gamma'(x^{m-1})}e^{-\lambda t_{k^*}}p_{k^*}\leq e^{-\lambda t_{k^*}}p_{k^*}.
\]

If $(i^*,k^*)\in\n T_0$, then $i^*=0, u_{0\bar t}(k^*)\geq0$, $u_{0x}(k^*)\leq0$. Notice that $\left(\frac1{\gamma_0}\right)_x = -\frac1{\gamma_0\gamma_{1}}\gamma_{0x}$. Note $\gamma_{0x}=\gamma'(x^0)$, so for $h$ small enough, we can ascertain $\gamma_{0x}=\gamma'(x^0)\in(\frac12,\frac32)$. It follows
\[
\frac{\gamma'(x^0)}{\gamma_0\gamma_1}u_0(k^*) \leq e^{-\lambda t_{k^*}}\big(hf_{0k^*}-g_{k^*}^n\big)\quad\implies\quad u_0(k^*)\leq e^{-\lambda t_{k^*}}\big(hf_{0k^*}-g_{k^*}^n\big).
\]

If $(i^*,k^*)\in\n N$, then $u_{i^*\bar t}(k^*)\geq 0$, $~u_{i^*x\bar x}(k^*)=\frac1{h^2} \big(u_{i^*+1}(k^*)-2u_{i^*}(k^*)+u_{i^*-1}(k^*)\big)\leq 0$. For $(i,k)\in\n N$, the corresponding equation in (\ref{systemu}) is equivalent to
\begin{gather}
\frac1{\gamma_i}\zeta_k^ie^{-\lambda\tau}u_{i\bar t}(k)-\frac1{\gamma_{i+1}}u_{ix\bar x}(k) -\left[ \left(\frac1{\gamma_i}\right)_{\bar x} + \left(\frac1{\gamma_i}\right)_{x} \right]u_{i\bar x}(k) ~\qquad\qquad\qquad\qquad \nonumber\\[3mm] \qquad\qquad\qquad +~ \left[\frac{\lambda}{\gamma_i}\zeta_k^ie^{-\lambda(t_k-t^k)} - \left(\frac1{\gamma_i}\right)_{x\bar x}\right]u_i(k)=f_{ik}e^{-\lambda t_k}. \label{other}
\end{gather}

Define the sets
\[
\n N_+=\left\{(i,k)\in\n N\Big| \left(\frac1{\gamma_i}\right)_{\bar x} + \left(\frac1{\gamma_i}\right)_{x} \geq 0 \right\},\quad \n N_-=\left\{(i,k)\in\n N\Big| \left(\frac1{\gamma_i}\right)_{\bar x} + \left(\frac1{\gamma_i}\right)_{x} < 0 \right\}.
\]
And it's clear $\n N=\n N_+\cup\n N_-$. Suppose $(i^*,k^*)\in\n N_+$. Then owing to (\ref{systemu}) and $u_{i^*x}(k^*)\leq0$, we can write
\begin{equation}\label{est}
\left[\frac{\lambda}{\gamma_{i^*}}\zeta_{k^*}^{i^*}e^{-\lambda(t_{k^*}-t^{k^*})} - \left(\frac1{\gamma_{i^*}}\right)_{x\bar x}\right]u_{i^*}(k^*)\leq f_{i^*k^*}e^{-\lambda t_{k^*}}.
\end{equation}
If instead $(i^*,k^*)\in\n N_-$, then we can use (\ref{other}) and the fact that $~u_{i^*\bar x}(k^*)\geq0$ to achieve again (\ref{est}). Therefore, (\ref{est}) is achieved in any case. We can choose $\tau$ so small that $e^{-\lambda(t_{k}-t^{k})}>\frac12,\quad\forall k$. Observe that
\begin{align*}
\frac1{\zeta_{k^*}^{i^*}}e^{\lambda(t_{k^*}-t^{k^*})}\left[-\frac1{\gamma_{i^*+1}}\gamma_{i^*x\bar x} + \frac{\gamma_{i^*x}+\gamma_{i^*\bar x}}{\gamma_{i^*-1}\gamma_{i^*+1}}\gamma_{i^*\bar x}\right] &\leq \frac2{\bar b}\left[4(2\Vert\gamma''\Vert_{C[0,\ell]})+ 16(2\Vert\gamma'\Vert_{C[0,\ell]}^2)\right] \nonumber \\ &\leq \frac{64}{\bar b} \big(\Vert\gamma''\Vert_{C[0,\ell]} +\Vert\gamma'\Vert_{C[0,\ell]}^2\big).
\end{align*}
Then by (\ref{lambda}), it is the case that the coefficient of $u_{i^*}(k^*)$ is positive independently of $i^*,k^*$. Therefore,
\[
u_{i^*}(k^*)\leq C_{\gamma}f_{i^*k^*}e^{-\lambda t_{k^*}}
\]
where $C_{\gamma}$ is a constant depending only on $\gamma$ and $\bar b$.

We can put together the obtained estimations to deduce that for $(i,k)\in \n M_{k_1}$,
\[
u_i(k)\leq\max\limits_{\n M}u_i(k) \leq\max\Big\{0,~~\Vert\Phi\Vert_{C[0,\ell]}, ~~\Vert p\Vert_{C[0,T]}, ~~\Vert g^n\Vert_{C[0,T]}+\Vert f\Vert_{L_{\infty}(D)}, ~~C_{\gamma}\Vert f\Vert_{L_{\infty}(D)}  \Big\}.
\]
But because $u_i(k) =\gamma_ie^{-\lambda t_k}v_i(k)$, we have the following uniform upper bound for the discrete state vector:
\[
v_i(k)\leq4e^{\lambda T}\max\Big\{0,~~\Vert\Phi\Vert_{C[0,\ell]}, ~~\Vert p\Vert_{C[0,T]}, ~~\Vert g^n\Vert_{C[0,T]}+\Vert f\Vert_{L_{\infty}(D)}, ~~C_{\gamma}\Vert f\Vert_{L_{\infty}(D)}  \Big\},\quad (i,k)\in \n M_{k_1}.
\]
In a fully analogous manner, we arrive at a uniform lower bound for the discrete state vector: For $(i,k)\in \n M_{k_1}$,
\[
v_i(k)\geq4e^{\lambda T}\min\Big\{0,~~-\Vert\Phi\Vert_{C[0,\ell]}, ~~-\Vert p\Vert_{C[0,T]}, ~~-\Vert g^n\Vert_{C[0,T]}-\Vert f\Vert_{L_{\infty}(D)}, ~~-C_{\gamma}\Vert f\Vert_{L_{\infty}(D)}  \Big\}.
\]
Combining the uniform upper and lower bounds imply (\ref{linfestimate}) up to $k_1$. But $k_1$ was arbitrary in $1,\ldots,n$. Theorem is proved. \hfill{$\square$}

\subsection{$W_{2}^{1,1}$- energy estimation for the Discrete Stefan Problem}\label{E:3:1b}

\begin{theorem}\label{energyestimate} Suppose that $p\in W_2^1(0,T), \Phi\in W_2^1(0,\ell), f\in L_{\infty}(D)$. For $[g]_n\in \n G_R^n$ and $n,m$ large enough, the discrete state vector $[v([g]_n)]_n$ satisfies the following estimate:
	\begin{gather}
	\label{energy}\Vert[v]_n\Vert^2_{\n E}:= \sum\limits_{k=1}^n\tau\sum\limits_{i=0}^{m-1}hv_{i\bar t}^2(k) +\max\limits_{1\leq k\leq n}\left(\sum\limits_{i=0}^{m-1}hv_{ix}^2(k)\right) + \sum\limits_{k=1}^n\tau^2 \sum\limits_{i=0}^{m-1} hv_{ix\bar t}^2(k)\\ \leq~~\tilde C_{\infty}\Big(\Vert\Phi\Vert^2_{W_2^1(0,\ell)} +\Vert f\Vert_{L_{\infty}(D)}^2 +\Vert p\Vert_{W_2^1(0,T)}^2 +\Vert g^n\Vert_{W_2^1(0,T)}^2\Big) \nonumber
	\end{gather}
	where $\tilde C_{\infty}$ is a constant independent of $n$ and $m$.
\end{theorem}

\noindent\emph{Proof. }Consider $n$ and $m$ large enough that Theorem \ref{uniformboundedness} is satisfied. In (\ref{dsvsum}), choose $\eta = 2\tau v_{\bar t}(k)$. Using (\ref{zetav}), write $(b_n(v_i(k)))_{\bar t} =\zeta_k^iv_{i\bar t}(k)$. Also, use the fact that
\begin{align*}
2\tau v_{ix}(k)(v_{i\bar t}(k))_x &= 2\tau v_{ix}(k)(v_{ix}(k))_{\bar t} \\[3mm] &= v_{ix}^2(k)+ v_{ix}^2(k)-2v_{ix}(k)v_{ix}(k-1)+v_{ix}^2(k-1)-v_{ix}^2(k-1) \\[3mm] &= v_{ix}^2(k)+\Big(v_{ix}(k)-v_{ix}(k-1)\Big)^2-v_{ix}^2(k-1) \\[3mm] &= v_{ix}^2(k)-v_{ix}^2(k-1) +\tau^2 v_{ix\bar t}^2(k).
\end{align*}
We thus have
\begin{gather}
2\tau\sum\limits_{i=0}^{m-1}h\zeta_k^iv_{i\bar t}^2(k) +\sum\limits_{i=0}^{m-1}hv_{ix}^2(k)-\sum\limits_{i=0}^{m-1}hv_{ix}^2(k-1)+ \tau^2\sum\limits_{i=0}^{m-1}hv_{ix\bar t}^2(k)\nonumber \\ = 2\tau\sum\limits_{i=0}^{m-1}hf_{ik}v_{i\bar t}(k)+2\tau p_kv_{m\bar t}(k)-2\tau g_k^nv_{0\bar t}(k) \label{e1}.
\end{gather}
Estimate the right-hand side of (\ref{e1}) by applying Cauchy Inequality with $\ep >0$ in the first term. Recall that $b_n'(v)\geq \bar b, ~\forall v$. We have:
\begin{gather}
2\tau\sum\limits_{i=0}^{m-1}h\zeta_k^iv_{i\bar t}^2(k) +\sum\limits_{i=0}^{m-1}hv_{ix}^2(k)-\sum\limits_{i=0}^{m-1}hv_{ix}^2(k-1)+ \tau^2\sum\limits_{i=0}^{m-1}hv_{ix\bar t}^2(k) \nonumber \\ \leq \bar b\tau\sum\limits_{i=0}^{m-1}hv_{i\bar t}^2(k) + \frac1{\bar b}\tau\sum\limits_{i=0}^{m-1}hf_{ik}^2 +2\tau p_kv_{m\bar t}(k)-2\tau g_k^nv_{0\bar t}(k) \label{e2}.
\end{gather}
We can absorb the first term on the right-hand side of (\ref{e2}) to the first term on the left-hand side. Hence:
\begin{gather}
\tau\sum\limits_{i=0}^{m-1}h\bar bv_{i\bar t}^2(k) +\sum\limits_{i=0}^{m-1}hv_{ix}^2(k)-\sum\limits_{i=0}^{m-1}hv_{ix}^2(k-1)+ \tau^2\sum\limits_{i=0}^{m-1}hv_{ix\bar t}^2(k) \nonumber \\ \leq \frac1{\bar b}\tau\sum\limits_{i=0}^{m-1}hf_{ik}^2 +2\tau p_kv_{m\bar t}(k)-2\tau g_k^nv_{0\bar t}(k),\qquad \forall k =1,\ldots,n \label{e3}.
\end{gather}
Perform summation of (\ref{e3}) for $k$ from $1$ to $q,~~2\leq q\leq n$. The second and third terms on the left-hand side telescope, and we obtain:
\begin{gather}
\bar b\sum\limits_{k=1}^q\tau\sum\limits_{i=0}^{m-1}hv_{i\bar t}^2(k) +\sum\limits_{i=0}^{m-1}hv_{ix}^2(q)+ \sum\limits_{k=1}^q\tau^2 \sum\limits_{i=0}^{m-1} hv_{ix\bar t}^2(k) \nonumber \\ \leq \sum\limits_{i=0}^{m-1}hv_{ix}^2(0) + \frac1{\bar b}\sum\limits_{k=1}^q\tau\sum\limits_{i=0}^{m-1}hf_{ik}^2 +2\sum\limits_{k=1}^q\tau p_kv_{m\bar t}(k)-2\sum\limits_{k=1}^q\tau g_k^nv_{0\bar t}(k) \label{aftertelescope}.
\end{gather}
Use the summation by parts technique on the $p$ and $g$ sums:
\begin{gather}
\sum\limits_{k=1}^q\tau p_kv_{m\bar t}(k) = \sum\limits_{k=1}^qp_kv_m(k)- \sum\limits_{k=1}^qp_kv_m(k-1) = \sum\limits_{k=1}^qp_kv_m(k)- \sum\limits_{k=0}^{q-1}p_{k+1}v_m(k)\nonumber \\ = -\sum\limits_{k=1}^{q-1}\tau p_{kt}v_m(k) + p_qv_m(q) - p_1v_m(0), \nonumber \\ 
\sum\limits_{k=1}^q\tau g_k^nv_{0\bar t}(k) = -\sum\limits_{k=1}^{q-1}\tau g_{kt}^nv_0(k) + g_q^nv_0(q) - g_1^nv_0(0) \label{sumparts}.
\end{gather}
In view of (\ref{sumparts}) and borrowing (\ref{linfestimate}) from Theorem \ref{uniformboundedness}, (\ref{aftertelescope}) yields (through Cauchy Inequality):
\begin{gather}
\bar b\sum\limits_{k=1}^q\tau\sum\limits_{i=0}^{m-1}hv_{i\bar t}^2(k) +\sum\limits_{i=0}^{m-1}hv_{ix}^2(q)+ \sum\limits_{k=1}^q\tau^2 \sum\limits_{i=0}^{m-1} hv_{ix\bar t}^2(k) \leq \nonumber \\ \leq \sum\limits_{i=0}^{m-1}h\Phi_{ix}^2 + \frac1{\bar b}\sum\limits_{k=1}^q\tau\sum\limits_{i=0}^{m-1}hf_{ik}^2 +\sum\limits_{k=1}^{q-1}\tau p_{kt}^2 + \sum\limits_{k=1}^{q-1}\tau \big(g_{kt}^n\big)^2 ~~+ \nonumber \\[3mm]\qquad +~~ 2t_{q-1}\Vert[v]_n \Vert_{ \ell_{\infty}}^2 +2 \big(\Vert p\Vert_{L_{\infty}(0,T)}+\Vert g^n\Vert_{L_{\infty}(0,T)}\big) \Vert[v]_n\Vert_{\ell_{\infty}} \label{e4}.
\end{gather}
Through the definition of the Steklov average, the Cauchy-Schwarz inequality and Fubini's Theorem, for $h$ small enough we have the following results:
\begin{gather}
\sum\limits_{i=0}^{m-1}h\Phi_{ix}^2 \leq \Vert\Phi'\Vert^2_{L_2(0,\ell)} + \Vert\Phi'\Vert_{L_2(\ell-h,\ell)}^2 \leq \Vert\Phi\Vert^2_{W_2^1(0,\ell)}, \nonumber \\
\Vert g^n\Vert_{L_{\infty}(0,T)}\leq \Vert g^n\Vert_{W_2^1(0,T)}\leq R, \nonumber \\
\sum\limits_{k=1}^{q-1}\tau \big(g_{kt}^n\big)^2 =\sum\limits_{k=2}^q\tau(g_{k\bar t}^n)^2 \leq \Vert (g^n)'\Vert_{L_2(0,T)}^2\leq \Vert g^n\Vert_{W_2^1(0,T)}^2\leq R^2, \nonumber \\ 
\sum\limits_{k=1}^{q-1}\tau p_{kt}^2 =\sum\limits_{k=2}^q\tau p_{k\bar t}^2 \leq \Vert p' \Vert_{L_2(0,T)}^2\leq \Vert p\Vert_{W_2^1(0,T)}^2, \nonumber \\ 
\sum\limits_{k=1}^q\tau\sum\limits_{i=0}^{m-1}hf_{ik}^2 \leq \Vert f\Vert_{L_2(D)}^2 \label{estimations}.
\end{gather}
Applying the results in (\ref{estimations}) to (\ref{e4}),
\begin{gather}
\sum\limits_{k=1}^q\tau\sum\limits_{i=0}^{m-1}hv_{i\bar t}^2(k) +\sum\limits_{i=0}^{m-1}hv_{ix}^2(q)+ \sum\limits_{k=1}^q\tau^2 \sum\limits_{i=0}^{m-1} hv_{ix\bar t}^2(k) \nonumber \\\leq \tilde C_{\infty}\Big(\Vert\Phi\Vert^2_{W_2^1(0,\ell)} +\Vert f\Vert_{L_\infty(D)}^2 +\Vert p\Vert_{W_2^1(0,T)}^2 +\Vert g^n\Vert_{W_2^1(0,T)}^2\Big) \label{e5}
\end{gather}
where $\tilde C_{\infty}$ is a constant dependent on $\bar b, T, R$, but independent of $n,m$ and $q$. Since, $q=\overline{1,n}$ is arbitrary, from (\ref{e5}), (\ref{energy}) follows.\hfill{$\square$}

\begin{theorem}\label{compactness} Let $\{[g]_n\}$ be a sequence in $\n G_R^n$ such that the sequence of interpolations $\{\n P_n([g]_n)\}$ converges weakly to $g\in W_2^1[0,T]$. Then the whole sequence of interpolations $\{\hat v^{\tau}\}$ of the associated discrete state vectors converges weakly to $v=v(x,t;g)\in W_2^{1,1}(D)$, with $v$ the unique weak solution to the Stefan Problem in the sense of (\ref{weaksol}).
\end{theorem}

\noindent\emph{Proof. } By the definitions of the interpolations given in (\ref{interpolations}), and by using (\ref{system}) we deduce that
\begin{gather}
\Vert \hat v^{\tau}\Vert_{L_{\infty}(D)} = \underset{(x,t)\in D}{\text{esssup}}~|\hat v^{\tau}(x,t)| = \max\limits_{0\leq k\leq n}\Big(\max\limits_{0\leq i\leq m}|v_i(k)|\Big) =	\Vert[v]_n\Vert_{\ell_{\infty}} \label{compactnesslinf2} ,\\[4mm]
\int\limits_0^T\int\limits_0^{\ell}( \hat v^{\tau})^2\,dx\,dt \leq T\ell\Vert \hat v^{\tau}\Vert^2_{L_{\infty}(D)} = T\ell\Vert[v]_n\Vert^2_{\ell_{\infty}},  \nonumber \\
\int\limits_0^T\int\limits_0^{\ell}(\hat v^{\tau}_x)^2\,dx\,dt \leq2\sum\limits_{k=0}^{n-1}\sum\limits_{i=0}^{m-1}\tau hv_{ix}^2(k) + 2\sum\limits_{k=1}^n\sum\limits_{i=0}^{m-1}\frac13\tau^3 hv_{ix\bar t}^2(k) \leq 2(\Vert[v]_n\Vert_{\n E}^2 +\Vert\Phi\Vert_{W_2^1(0,\ell)}^2),\nonumber \\
\int\limits_0^T\int\limits_0^{\ell}(\hat v^{\tau}_t)^2\,dx\,dt \leq 2\sum\limits_{k=1}^n\sum\limits_{i=0}^{m-1}\tau h\big(v_{i\bar t}^2(k)+\frac13h^2v_{ix\bar t}^2(k)\big) \leq \frac23\sum\limits_{k=1}^n\left[\sum\limits_{i=0}^{m-1}\Big(7\tau hv_{i\bar t}^2(k)\Big)+2\tau hv_{m\bar t}^2(k)\right] \nonumber, \\
\sum\limits_{k=1}^n\tau hv_{m\bar t}^2(k) \leq 2\sum\limits_{k=1}^n\tau h\Big(v_{m-1,\bar t}^2(k)+h^2p_{k\bar t}^2\Big) \leq 2\Vert[v]_n\Vert_{\n E}^2+2h^3\Vert p\Vert_{W_2^1(0,T)}^2 \label{compactnessw211}.
\end{gather}
Since $[g]_n\in\n G^n$, then $\Vert g^n\Vert_{W_2^1[0,T]}\leq R+1$ for large enough $n$. From the energy estimates (\ref{linfestimate}), (\ref{energy}) and calculations (\ref{compactnesslinf2}), (\ref{compactnessw211})  it is therefore the case that $\{\hat v^{\tau}\}$ is uniformly bounded in the spaces $W_2^{1,1}(D)$ and $L_{\infty}(D)$. As such, we may choose a subsequence of $\{\hat v^{\tau}\}$ that converges weakly in $W_2^{1,1}(D)$ to some function $v\in W_2^{1,1}(D)\cap L_{\infty}(D)$, and thus strongly in $L_2(D)$, by virtue of which we can choose a further subsequence that converges to $v$ pointwise almost everywhere. It is our intent to show now that $v$ satisfies the definition of a weak solution to the Stefan Problem. To do this, first realize that the sequences $\{v^{\tau}\},\{\hat v^{\tau}\}$ are equivalent in $W_2^{1,0}(D)$, and sequences $\{v^{\tau}\},\{\tilde v\}$ are equivalent in $L_2(D)$, as shown by the following calculations:
\begin{gather}
\int\limits_0^T\int\limits_0^{\ell}|v^{\tau}-\hat v^{\tau}|^2\,dx\,dt = \sum\limits_{k=1}^n\sum\limits_{i=0}^{m-1}\int\limits_{t_{k-1}}^{t_k} \int\limits_{x_i}^{x_{i+1}}\big(\hat v(x;k)-\hat v(x;k-1)-\hat v_{\bar t}(x;k)(t-t_{k-1})\big)^2\,dx\,dt\nonumber \\ = \sum\limits_{k=1}^n \sum\limits_{i=0}^{m-1} \int\limits_{t_{k-1}}^{t_k} \int\limits_{x_i}^{x_{i+1}}\hat v_{\bar t}^2(x;k)(t_k-t)^2\,dx\,dt = \frac13\tau^2 \int\limits_0^T\int\limits_0^{\ell}(\hat v^{\tau}_t)^2\,dx\,dt\longrightarrow 0 \label{equivl2}, \\
\int\limits_0^T\int\limits_0^{\ell}|v^{\tau}_x-\hat v^{\tau}_x|^2\,dx\,dt = \sum\limits_{k=1}^n\sum\limits_{i=0}^{m-1}\int\limits_{t_{k-1}}^{t_k} \int\limits_{x_i}^{x_{i+1}}\big(v_{ix}(k)-v_{ix}(k-1)-v_{ix\bar t}(k)(t-t_{k-1})\big)^2\,dx\,dt\nonumber \\ = \sum\limits_{k=1}^n \sum\limits_{i=0}^{m-1} \int\limits_{t_{k-1}}^{t_k}hv_{ix\bar t}^2(k)(t_k-t)^2\,dt = \frac13\tau\left( \sum\limits_{k=1}^n\tau^2\sum\limits_{i=0}^{m-1}hv_{ix\bar t}^2(k)\right)\longrightarrow 0 \label{equivw210},
\end{gather}
\begin{gather}
\int\limits_0^T\int\limits_0^{\ell}|v^{\tau}-\tilde v|^2\,dx\,dt =\sum\limits_{k=1}^n\sum\limits_{i=0}^{m-1}\tau\int\limits_{x_i}^{x_{i+1}}\Big|v_i(k)+v_{ix}(k)(x-x_i)-v_i(k)\Big|^2\,dx\nonumber \\ = \sum\limits_{k=1}^n\sum\limits_{i=0}^{m-1} \frac13\tau h^3 v_{ix}^2(k) \longrightarrow 0 \label{equivconst},
\end{gather}
as $n,m$ go to $+\infty$. Accordingly, $v^{\tau}\rightarrow v$ weakly in $W_2^{1,0}(D)$ and $\tilde v\rightarrow v$ strongly in $L_2(D)$ and pointwise a.e. on $D$ along a subsequence. Fix arbitrary $\psi\in W_2^{1,1}(D)$ with $\psi|_{t=T}=0$. Actually, due to density of $C^1(\bar{D})$ in $W_2^{1,1}(D)$, without loss of generality we can consider $\psi\in C^1(\overline D)$ and $\psi|_{t=T}=0$. Define $\psi_i(k) = \psi(x_i,t_k),~\forall i~\forall k$, and consider the interpolations:
\begin{gather}
\psi^{\tau}(x,t):= \psi_i(k),\qquad \psi_x^{\tau}(x,t):= \psi_{ix}(k)\qquad \psi_t^{\tau}(x,t):= \psi_{it}(k), \nonumber \\ \qquad \qquad x_i\leq x<x_{i+1},\quad t_{k-1}< t\leq t_{k},\qquad i =\overline{0,m},~k =\overline{0,n} \label{psi}.
\end{gather}
It is readily checked that $\psi^{\tau},\psi_x^{\tau},\psi_t^{\tau}$ converge uniformly on $\overline D$ as $n,m\rightarrow\infty$ to the functions $\psi, \psi_x,\psi_t$ respectively. Fix $n$. For each $k$ in (\ref{dsvsum}) as satisfied by the discrete state vector $[v([g]_n)]_n$, choose $\eta_i=\tau\psi_i(k),~\forall i$ and sum all equalities (\ref{dsvsum}) over $k=1,\ldots,n$. The resulting expression is as follows:
\begin{gather}
\sum\limits_{k=1}^n\tau\sum\limits_{i=0}^{m-1}h\Big[\big(b_n(v_i(k))\big)_{\bar t}\psi_i(k) + v_{ix}(k)\psi_{ix}(k)-f_{ik}\psi_i(k)\Big] -\sum\limits_{k=1}^n\tau p_k\psi_m(k) +\sum\limits_{k=1}^n\tau g_k^n\psi_0(k). \label{c1}
\end{gather}
We transform the first term through summation by parts:
\begin{gather*}
\sum\limits_{k=1}^n\tau\sum\limits_{i=0}^{m-1}h\big(b_n(v_i(k))\big)_{\bar t}\psi_i(k) =\sum\limits_{k=1}^n\sum\limits_{i=0}^{m-1}hb_n(v_i(k))\psi_i(k)-\sum\limits_{k=1}^n\sum\limits_{i=0}^{m-1}hb_n(v_i(k-1))\psi_i(k) \nonumber \\ = \sum\limits_{k=1}^n\sum\limits_{i=0}^{m-1}hb_n(v_i(k))\psi_i(k)-\sum\limits_{k=0}^{n-1}\sum\limits_{i=0}^{m-1}hb_n(v_i(k))\psi_i(k+1) \nonumber \\ = -\sum\limits_{k=1}^{n-1}\tau \sum\limits_{i=0}^{m-1}hb_n(v_i(k))\psi_{it}(k) + \sum\limits_{i=0}^{m-1}hb_n(v_i(n))\psi_i(n) - \sum\limits_{i=0}^{m-1}hb_n(v_i(0))\psi_i(1)  \nonumber \\ = -\int\limits_0^{T-\tau}\int\limits_0^{\ell}b_n(\tilde v(x,t))\psi_t^{\tau}\,dx\,dt -\int\limits_0^{\ell}b_n(\tilde\Phi(x))\psi^{\tau}(x,\tau)\,dx.
\end{gather*}
Thus, (\ref{c1}) can be rewritten as:
\begin{gather}
\int\limits_0^T\int\limits_0^{\ell}\Big[-b_n(\tilde v)\psi_t^{\tau} +v_x^{\tau}\psi_x^{\tau} -f\psi^{\tau}\Big]\,dx\,dt -\int\limits_0^{\ell}b_n(\tilde\Phi)\psi^{\tau}(x,\tau)\,dx  \nonumber \\ -\int\limits_0^Tp(t)\psi^{\tau}(\ell,t)\,dt +\int\limits_0^Tg^n(t)\psi^{\tau}(0,t)\,dt + \int\limits_{T-\tau}^T\int\limits_0^{\ell}b_n(\tilde v)\psi_t^{\tau}\,dx\,dt = 0 \label{c2}.
\end{gather}
Theorem \ref{uniformboundedness} implies that if $\n V_n:=\big\{y\in\bb R~|~\exists(x,t)\in D~\text{s.t. }\tilde v(x)=y\}$ (i.e. $\n V_n$ is the range of $\tilde v$), then the set $\n V = \bigcup\limits_{n=1}^{\infty}\n V_n$ is bounded in $\bb R$, hence its closure $\overline{\n V}$ is compact in $\bb R$. Because of the piecewise continuity of $b$, it follows that $b(\tilde v(x,t))\in L_{\infty}(D)$, and therefore $\Vert b_n(\tilde v)\Vert_{L_{\infty}(D)} \leq C$. Since $D$ is a set of finite measure, $\Vert b_n(\tilde v)\Vert_{L_{2}(D)} \leq C$, so that a subsequence $\{b_{n_l}(\tilde v(x,t))\}$ can be constructed so that it converges weakly in $L_2(D)$ to a function $\tilde b(x,t)\in L_2(D)$. Through a similar argument, we can choose this subsequence so that $b_{n_l}(\tilde\Phi(x))$ converges weakly in $L_2[0,\ell]$ to a function $\tilde b_0(x)\in L_2[0,\ell]$. Take a diagonal of these subsequences as the whole sequence. We see that
\begin{align}
\int\limits_{T-\tau}^T\int\limits_0^{\ell}b_n(\tilde v)\psi_t^{\tau}\,dx\,dt &\leq \left(~\int\limits_{T-\tau}^T\int\limits_0^{\ell}b_n^2(\tilde v)\,dx\,dt\right)^{\frac12} \left(~\int\limits_{T-\tau}^T\int\limits_0^{\ell}(\psi_t^{\tau})^2\,dx\,dt\right)^{\frac12} \nonumber \\[4mm] &\leq \Vert b_n^2(\tilde v)\Vert_{L_2(D)} \Vert\psi_t^{\tau}\Vert_{L_2\big([0,\ell]\times[T-\tau,T]\big)}\longrightarrow0~\text{as } n\rightarrow\infty\label{ex1}.
\end{align}
Now, due to (\ref{ex1}), the uniform convergence of $\psi^{\tau},~\psi_x^{\tau},~\psi_t^{\tau}$ respectively to $\psi, ~\psi_x,~\psi_t$ and weak convergence of $b_n(\tilde v), ~v_x^{\tau},$ $~b_n(\tilde\Phi), ~g^n$ to $\tilde b,~ v_x,~ \tilde b_0,~ g$ in the respective $L_2$ spaces, then as $n\rightarrow\infty$, (\ref{c2}) implies
\begin{gather}
\int\limits_0^T\int\limits_0^{\ell}\Big[-\tilde b(x,t)\psi_t +v_x\psi_x -f\psi\Big]\,dx\,dt -\int\limits_0^{\ell}\tilde b_0(x)\psi(x,0)\,dx  \nonumber \\ -\int\limits_0^Tp(t)\psi(\ell,t)\,dt +\int\limits_0^Tg(t)\psi(0,t)\,dt = 0 \label{vlimweaksol}.
\end{gather}
It can be checked that both $\tilde b$ and $\tilde b_0$ are functions of type $\n B$. If at the point $(x,t)$,
\[ \tilde{v}(x,t) \rightarrow v(x,t) \neq v^j  ,\]
then we have
\[
b_n(\tilde v(x,t)) = \int\limits_{\tilde v(x,t)-\frac1n}^{\tilde v(x,t)+\frac1n}\omega_{1/n}(|\tilde v(x,t)-u|)b(u)\,du~\longrightarrow~b(v(x,t)).
\]
On the contrary, if at the point $(x,t)$ we have
\[ \tilde{v}(x,t) \rightarrow v(x,t) = v^j , \]
then we have
\[ b(v^j)^- \leq \liminf_{n\rightarrow \infty}b_n(\tilde v(x,t)) \leq \limsup_{n\rightarrow \infty}b_n(\tilde v(x,t)) \leq b(v^j)^+ . 
\]
Since the sequence $\{b_n(\tilde{v})\}$ converges to $\tilde{b}(x,t)$ weakly in $L_2(D)$, by Mazur's lemma there is a sequence of convex combinations of 
elements of $\{b_n(\tilde{v})\}$ which converges to $\tilde{b}(x,t)$ strongly in $L_2(D)$. Therefore, there is a subsequence of convex combinations which converges to $\tilde{b}(x,t)$ a.e. in $D$. It easily follows that  $\tilde b=B(x,t,v(x,t))$ is a function of type $\n B$. In a very similar way, it is seen that $\tilde b_0=B(x,0,\Phi(x))$ is of type $\n B$. Hence, by definition, $v$ is a weak solution to the Stefan Problem in the sense of (\ref{weaksol}). From Lemma \ref{unique} then, $v$ is the unique solution, which implies that $v$ is the only weak limit point of the sequence $\{\hat v^{\tau}\}$. Therefore, the whole sequence $\{\hat v^{\tau}\}$ converges to $v$ weakly in $W_2^{1,1}(D)$. \hfill{$\square$}

\subsection{Existence of the Optimal Control}\label{E:3:1c}
Consider a sequence $\{g_l\}\in \n G_R$ such that $\n J(g_l)\searrow \n J_*$. Since $\{g_l\}$ is uniformly bounded in $W_2^1(0,T)$, it is weakly precompact in $\n G_R$. Therefore, there exists a subsequence $\{g_{l_k}\}$ which converges weakly in $W_2^1(0,T)$, say, to $g\in W_2^1(0,T)\in \n G_R$. For ease of notation, take this subsequence as the sequence $\{g_l\}$.

Let $v_l = v(x,t;g_l)$ and $v=v(x,t;g)$ be solutions to the Stefan problem in the sense of (\ref{weaksol}) with $g_l$ and $g$ respectively. Then for fixed $l$, the sequence of vectors $\{[g_l]_n\}$ given by $[g_l]_n =\n Q_n(g_l)$ is such that the interpolations $g_l^n =\n P_n([g_l]_n)$ converge weakly in $W_2^1(0,T)$ to $g_l\in W_2^1(0,T)$ as $n\rightarrow\infty$. Therefore, Theorem \ref{compactness} applies, and so associated to $[g_l]_n$ the interpolations $\hat v^{\tau}_l$ of the discrete state vectors $[v([g_l]_n)]_n$ converge weakly in $W_2^{1,1}(D)$ to $v_l$. As such,
\begin{equation}\label{weakcomp}
\Vert v_l\Vert_{W_2^{1,1}(D)} \leq \liminf\limits_{n\rightarrow\infty}\Vert\hat v^{\tau}_l\Vert_{W_2^{1,1}(D)} \leq C\liminf\limits_{n\rightarrow\infty} \Big(\Vert[v_l]_n\Vert_{\ell_{\infty}} + \Vert[v_l]_n\Vert_{\n E}\Big)
\end{equation}
where $C$ is independent of $n, m$ and $l$. Thanks to $\{g_l\}\subset\n G_R$, it it is clear from (\ref{linfestimate}) and (\ref{energy}) that the right-hand side of (\ref{weakcomp}) is uniformly bounded. Similarly, one can conclude that
\[
\Vert v_l\Vert_{L_{\infty}(D)} \leq\liminf\limits_{n\rightarrow\infty} \Vert[v_l]_n\Vert_{\ell_{\infty}}.
\]
Accordingly, $\{v_l\}\in W_2^{1,1}(D)\cap L_{\infty}(D)$ is a weakly precompact sequence in $W_2^{1,1}(D)$, so that it contains a subsequence $\{v_{l_k}\}$ which converges weakly to a function $\tilde v\in W_2^{1,1}(D)$, and thus strongly in $L_2(D)$. Due to this strong convergence in $L_2(D)$, a further subsequence of $\{v_{l_k}\}$ can be extracted which converges almost everywhere to $\tilde v$ on $D$. Then the uniform boundedness of this subsequence in $L_{\infty}(D)$ implies that $\tilde v\in L_{\infty}(D)$. Now, take this subsequence as the whole sequence. Each of the $v_l$ satisfies (\ref{weaksol}) with $g_l$ and with an arbitrarily fixed function $B$ of type $\n B$. Going to infinity along the sequence, we have that we can replace $g_l$ with $g$ and $v_l$ with $\tilde v$ in (\ref{weaksol}). Indeed, $B(x,t,v_l(x,t))\rightarrow B(x,t,\tilde v(x,t))$ a.e. on $D$ because of Corollary \ref{Cjmeasure0} and the fact that $v_l\rightarrow\tilde v$ a.e. on $D$. Consequently, $\tilde v$ is a solution to the Stefan problem with $g$. But, due to uniqueness of such a solution, it follows that $v=\tilde v$ in $W_2^{1,1}(D)\cap L_{\infty}(D).$ Next, note:
\begin{align}
\lim\limits_{l\rightarrow\infty}\Big|\n J(g)- \n J(g_l)\Big| &= \lim\limits_{l\rightarrow\infty}\Big| \Vert v(\ell,t)-\Gamma(t)\Vert ^2_{L_2[0,T]}-  \Vert v_l(\ell,t)-\Gamma(t)\Vert^2_{L_2[0,T]}\Big| \nonumber \\ &= \lim\limits_{l\rightarrow\infty}\Big|<v-\Gamma,v-\Gamma>_{L_2[0,T]} - <v_l-\Gamma,v_l-\Gamma>_{L_2[0,T]}\Big| \nonumber \\ &=\lim\limits_{l\rightarrow\infty}\Big|\Vert v(\ell,t)-v_l(\ell,t) \Vert_{L_2[0,T]}^2 + 2<v-v_l,v_l-\Gamma>_{L_2[0,T]}\Big| \nonumber \\ & = \lim\limits_{l\rightarrow\infty} \left|\int\limits_0^T|v(\ell,t)- v_l(\ell,t)|^2\,dt + 2\int\limits_0^T\big(v(\ell,t)-v_l(\ell,t)\big) \big(v_l(\ell,t)-\Gamma(t)\big)\,dt\right|\label{wcj}.
\end{align}
By the weak convergence of the sequence $\{v_l\}$ to $v$ in $W_2^{1,1}(D)$, it follows that we have strong convergence in the space of traces. In particular, the integrals in (\ref{wcj}) vanish as $l\rightarrow\infty$. Hence $\lim\limits_{l\rightarrow\infty}\n J(g_l) = \n J(g)$. This limit is unique though, therefore it is the case that $\n J(g) = \n J_*$, so that $g\in \n G_*$. \hfill{$\square$}

\subsection{Proof of the Convergence of Discrete Optimal Control Problem}\label{E:3:1d}
The proof of Theorem \ref{convergence} is split into three separate lemmas, as shown below.

\noindent\textbf{Lemma A.} \emph{Let $\n J_*(\pm\ep)=\inf\limits_{\n G_{R\pm\ep}}\n J(g),~\ep>0.\quad$ Then $\quad\lim\limits_{\ep\rightarrow0}\n J_*(\ep)=\n J_* = \lim\limits_{\ep\rightarrow0}\n J_*(-\ep)$.}

\noindent\emph{Proof. } The proof of this lemma is very similar to the analogous lemma from \cite{Abdulla1}. If $0<\ep_1<\ep_2$, then
\[\n J_*(\ep_2)\leq\n J_*(\ep_1)\leq\n J_*\leq \n J_*(-\ep_1)\leq \n J_*(-\ep_2)
\]
Hence $\lim\limits_{\ep\rightarrow0}\n J_*(\ep)\leq\n J_*$ and  $\lim\limits_{\ep\rightarrow0}\n J_*(-\ep)\geq\n J_*$. Choose $g_{\ep}\in\n G_{R+\ep}$ such that 
\[\lim\limits_{\ep\rightarrow0}\Big(\n J(g_{\ep})-\n J_*(\ep)\Big)= 0.
\]
Since $\{g_{\ep}\}$ is weakly pre-compact in $W_2^1[0,T]$, there exists a subsequence $\ep'$ such that $g_{\ep'}\rightarrow g_*$ weakly in $W_2^1[0,T]$ as $\ep'\rightarrow0$. Since $\n J$ is weakly continuous, $\n J(g_{\ep'})\rightarrow\n J(g_*)$ as $\ep\rightarrow0$. Hence $\n J_*(\ep')\rightarrow\n J(g_*)\geq\n J_*$ as $\ep'\rightarrow0$. Thus $\lim\limits_{\ep\rightarrow0}\n J_*(\ep)=\n J_*$.

From the other side, by Theorem \ref{existence} we know there exists $g_*\in\n G_R$ such that $\n J(g_*)=\n J_*$. If $g_*\in\n G_R\backslash\n\partial\n G_R$, then there exists $\ep_*>0$ such that $g_*\in\n G_{R-\ep},~\forall\ep<\ep_*$, and in this case $\n J_*(-\ep) =\n J_*,~\forall\ep<\ep_*$. If $g_*\in\partial\n G_R$, then there exists $\{g_{\ep}\}$ with $g_{\ep}\in\n G_{R-\ep}$ such that $g_{\ep}\rightarrow g_*$ in $W_2^1[0,T]$ as $\ep\rightarrow0$. The continuity of $\n J$ gives us that $\lim\limits_{\ep\rightarrow0}\n J(g_{\ep}) = \n J(g_*)=\n J_*$. Since on the other hand, $\n J(g_{\ep})\geq\n J_*(-\ep)$, it follows that $\lim\limits_{\ep\rightarrow0}\n J_*(-\ep) = \n J_*$.\hfill{$\square$}

\noindent\textbf{Lemma B.} \emph{For all $g\in\n G_R, \lim\limits_{n\rightarrow\infty}\n I_n(\n Q_n(g)) =\n J(g)$}.

\noindent\emph{Proof. } Take $g\in\n G_R$ arbitrarily. If $\n Q(g)=[g]_n$, and $g^n = \n P_n([g]_n)$, then $g^n\rightarrow g$ strongly in $W_2^1(0,T)$ as $n\rightarrow\infty$. Applying Theorem \ref{compactness}, we have that the interpolations $\hat v^{\tau}$ of the discrete state vectors $[v([g]_n)]_n$ converge to $v=v(x,t;g)$ weakly in $W_2^{1,1}(D)$ as $n\rightarrow\infty$, and thus the traces $\hat v^{\tau}(\ell,\cdot)$ converge strongly in $L_2(0,T)$ to trace $v(\ell,\cdot)$. By calculations (\ref{equivl2}) and (\ref{equivw210}), the sequences $\{v^{\tau}\}$, $\{\hat v^{\tau}\}$ are equivalent in $W_2^{1,0}(D)$, so that the $v^{\tau}(\ell,\cdot)$ traces too converge to $v(\ell,\cdot)$ strongly in $L_2(0,T)$. If we define
\begin{equation}\label{tildegamma}
\tilde\Gamma(t) = \Gamma_k=\frac1{\tau}\int\limits_{t_{k-1}}^{t_k}\Gamma(t)\,dt,\qquad t_{k-1}<t\leq t_k,\quad k=\overline{1,n},
\end{equation}
then $\tilde\Gamma\rightarrow\Gamma$ in $L_2(0,T)$ as $n\rightarrow\infty$. Therefore,
\begin{align}
|\n I_n(\n Q_n(g))-\n J(g)| &= \left|\sum\limits_{k=1}^n\tau(v_m(k)-\Gamma_k)^2~ -\int\limits_0^T(v(\ell,t)-\Gamma(t))^2\,dt\right| \nonumber \\ &\leq \Vert v^{\tau}(\ell,\cdot)-v(\ell,\cdot)\Vert_{L_2[0,T]}^2+\Vert\tilde\Gamma-\Gamma\Vert_{L_2[0,T]}^2\quad \nonumber \\ &~+ 2\int\limits_0^T\Big[|v^{\tau}(\ell,t)-v(\ell,t)| |v(\ell,t)-\tilde\Gamma(t)|+|\Gamma(t)-\tilde\Gamma(t)||v(\ell,t)-\Gamma(t)|\Big]\,dt \nonumber \\ &\qquad\longrightarrow0\quad\text{as}\quad n\rightarrow\infty,\nonumber
\end{align}
which establishes the lemma. \hfill{$\square$}

\noindent\textbf{Lemma C.} \emph{For an arbitrary sequence $\{[g]_n\}$ such that $[g]_n\in\n G_R^n,$
	\[
	\lim\limits_{n\rightarrow\infty}\Big(\n J(\n P_n([g]_n))-\n I_n([g]_n)\Big) = 0.
	\]}

\noindent\emph{Proof. } Let $g^n = \n P_n([g]_n)$. This sequence is weakly precompact, so that a subsequence $g^{n_l}$ converges weakly to a function $g$ in $W_2^1(0,T)$. Take this subsequence as the whole sequence. Note that
\begin{equation}
|\n J(g^n)-\n I_n([g]_n)|\leq |\n J(g^n)-\n J(g)|+|\n J(g)-\n I_n([g]_n)| \label{mainc}.
\end{equation}
Since $\n J$ is weakly continuous, $|\n J(g^n)-\n J(g)|\longrightarrow0$ as $n\rightarrow\infty$. It remains to show that the second term on the right-hand side of (\ref{mainc}) goes to $0$ as $n\rightarrow\infty$. Actually, the proof of this fact follows in a manner very similar to the proof of Lemma B. From (\ref{mainc}) it follows that $\lim\limits_{n_l\rightarrow\infty}\Big(\n J(\n P_{n_l}([g]_{n_l}))-\n I_{n_l}([g]_{n_l})\Big) = 0$. However, the subsequence chosen was arbitrary. Therefore, the same result can be achieved for any subsequence $\{g_{n_{\alpha}}\}$ of $\{g_n\}$. It is then the case that the whole sequence $\n J(\n P_n([g]_n))-\n I_n([g]_n)$ converges to $0$ as $n\rightarrow\infty$. \hfill{$\square$}


\end{document}